\UseRawInputEncoding   % required by arxiv
\documentclass{article}

\usepackage{amsthm,amsfonts,amssymb,amsmath}
\usepackage{enumerate}
\usepackage[colorlinks=true]{hyperref}

\newtheorem{thm}{Theorem}[section]
\newtheorem{cor}[thm]{Corollary}
\newtheorem{lem}[thm]{Lemma}

\newtheorem{pro}[thm]{Proposition}
\newtheorem{conj}[thm]{Conjecture}

\theoremstyle{remark}
\newtheorem{remark}[thm]{Remark}

\theoremstyle{definition}

\newtheorem{example}[thm]{Example}

\newcommand{\CC}{\mathbb{C}}
\newcommand{\RR}{\mathbb{R}}
\newcommand{\ZZ}{\mathbb{Z}}
\newcommand{\QQ}{\mathbb{Q}}

\newcommand{\PP}{\mathbb{P}}
\newcommand{\DD}{\mathbb{D}}

\newcommand{\CA}{{\mathcal {A}}}

\newcommand{\CF}{{\mathcal {F}}}
\newcommand{\CG}{{\mathcal {G}}}

\newcommand{\CL}{{\mathcal {L}}}

\newcommand{\CO}{{\mathcal {O}}}

\newcommand{\CS}{{\mathcal {S}}}
\newcommand{\CT}{{\mathcal {T}}}

\newcommand{\CX}{{\mathcal {X}}}
\newcommand{\CY}{{\mathcal {Y}}}

\newcommand{\Aut}{{\mathrm{Aut}}}

\newcommand{\Hom}{{\mathrm{Hom}}}

\renewcommand{\Im}{{\mathrm{Im}}}

\newcommand{\Gal}{{\mathrm{Gal}}}

\newcommand{\trdeg}{{\mathrm{trdeg}}}

\newcommand{\Pic}{\mathrm{Pic}}

\renewcommand{\Re}{{\mathrm{Re}}}
\newcommand{\reg}{{\mathrm{reg}}}

\newcommand{\sing}{{\mathrm{sing}}}

\newcommand{\Sp}{{\mathrm{Sp}}}

\newcommand{\Sym}{{\mathrm{Sym}}}

\newcommand{\tr}{{\mathrm{tr}}}

\DeclareMathOperator{\Spec}{Spec}

\renewcommand{\d}{\textnormal{d}}

\renewcommand{\i}{\mathrm{i}}

\newcommand{\wt}{\widetilde}

\newcommand{\pair}[1]{\langle {#1} \rangle}

\newcommand{\ds}{\displaystyle}

\newcommand{\ol }{\overline}

\newcommand{\lra}{\longrightarrow}

\newcommand{\ft}{\mathfrak{t}}
\newcommand{\fs}{\mathfrak{s}}

\newcommand{\OD}{{\overline{D}}}

\newcommand{\C}{{\mathbb C}}
\newcommand{\D}{{\mathbb D}}

\newcommand{\R}{{\mathbb R}}

\begin{document}
%------------------------------------------------------

%\title{The geometric Bombieri--Lang conjecture for varieties of maximal Albanese dimension}

\title{Partial Heights, Entire Curves, and the Geometric Bombieri--Lang Conjecture}

\author{Junyi Xie, Xinyi Yuan}

\maketitle
\tableofcontents

\section{Introduction}

The Mordell conjecture, proved by Faltings in 1983, is a milestone in the history of Diophantine geometry. The Bombieri--Lang conjecture, a high-dimensional generalization of the Mordell conjecture, asserts that a projective variety over a number field satisfying a reasonable hyperbolicity condition has only finitely many rational points. 
If we change the hyperbolicity condition to the condition of being general type, then the conjecture asserts that the set of rational points is not Zariski dense.

In Vojta's landmark Ph.D. thesis (cf. \cite{Vojta1987}), he formulated a highly non-trivial connection between Diophantine geometry and Nevanlinna theory, which leads to far-reaching conjectures. This connection is now called Vojta's dictionary. According to Vojta's dictionary,  the analogy of an infinite sequence of distinct rational points is an entire curve. This predicts the Bombieri--Lang conjecture.

Besides the Mordell conjecture and more generally the case for projective varieties which admit a finite morphism to the moduli stack of polarized abelian varieties proved by Faltings \cite{Faltings1983}, the only known case of the Bombieri--Lang conjecture over number fields is for subvarieties of abelian varieties. This was proved by Faltings \cite{Faltings1991, Faltings1994} by extending the proof of Vojta \cite{Vojta91} on the Mordell conjecture. 

The geometric Bombieri--Lang conjecture is an analogue of the Bombieri-Lang conjecture over function fields. For a philosophy, we
quote the following statement from Lang \cite[p. 781]{Lang74}:
\begin{quote}
\emph{As usual, the absolute Mordell property has a relative formulation for
algebraic families of hyperbolic varieties: If there is an infinity of sections,
then the family contains split subfamilies, and almost all sections are due
to constant sections.}
\end{quote}
A new phenomena is caused by split subfamilies (constant subvarieties).
A similar philosophy also applies to varieties of general type over function fields. 

The geometric Bombieri--Lang conjecture was previously proved in the following cases (in various versions). 
For curves, it was proved by Manin \cite{Manin63} (for characteristic 0), Grauert \cite{Grauert65} (for characteristic 0), and Samuel \cite{Samuel66} (for positive characteristics).
For subvarieties of abelian varieties, it was proved by Raynaud \cite{Raynaud83} (for characteristic 0), Buium \cite{Buium92} (for characteristic 0), Abramovich--Voloch \cite{AV92} (for positive characteristics), and Hrushovski \cite{Hrushovski96} (for all characteristics). 
For smooth projective varieties with ample cotangent bundles, it was proved by 
Noguchi \cite{Nog82} (for characteristic 0),
Martin-Deschamps \cite{MD84} (for characteristic 0),
and Gillet--R\"ossler \cite{GR18} (for all characteristics).
%\cite{Nog92} for trivial family, and for hyperbolically embedded case

In this paper, we introduce a new approach to the geometric Bombieri--Lang conjecture for hyperbolic varieties in characteristic 0. 
The idea is as follows.
Let $X$ be a projective variety over a function field $K$ of one variable over a field $k$ of characteristic $0$.  
By the Lefschetz principle, we can reduce the problem to the case $k=\CC$. 
Let $B$ be the complex projective curve with function field $K$, and $\CX$ be a (projective) integral model of $X$ over $B$.  
If $X(K)$ has a sequence of unbounded Weil heights, 
assuming a \emph{non-degeneracy conjecture}, 
we construct an \emph{entire curve} in a smooth closed fiber of $\CX$ over $B$, which implies that the fiber is not hyperbolic and thus induces a contradiction.
The construction relies on the \emph{classical Brody lemma} in complex geometry, which constructs an entire curve in a complex analytic variety as a limit of a suitable sequence of complex discs. Another ingredient of the approach is \emph{partial heights} of algebraic points on projective varieties over complex function fields, a new analytic notion generalizing the notion of Weil heights.
The non-degeneracy conjecture is a technical conjecture which asserts that if a sequence has unbounded Weil height, then it also has unbounded partial heights.

Our construction realizes the mechanism of Vojta's dictionary relating Diophantine geometry to Nevanlinna theory in a reasonably concrete way. 
Recall that a key part of the dictionary is an analogue between an infinite sequence of rational points and an entire curve.
To get a clean statement, we assume that $k$ is finitely generated over $\QQ$ in the above construction.
Then from any sequence of points in $X(K)$ with unbounded heights, our conditional
construction gives an entire curve in $X_\iota(\C)$ for ``most'' embeddings $\iota:K\hookrightarrow \C$. 

As an unconditional result, we prove the geometric Bombieri--Lang conjecture for any hyperbolic projective variety with a finite morphism to an abelian variety over a function field of characteristic $0$. 
In fact, we prove the non-degenerate conjecture in this case by applying Betti forms of complex abelian schemes introduced by Ngaiming Mok. 
This unconditional result includes the case of subvarieties of abelian varieties, and our approach is different from the previous works.

\subsection{The geometric Bombieri--Lang conjecture}

To state our main results, we will first formulate a precise statement of the geometric Bombieri--Lang conjecture. 
Although our main theorems are only for the hyperbolic case, we will formulate the conjecture in the general case.
We will start with some definitions related to hyperbolicity and Diophantine geometry. We refer to \cite{Lang87, Kob98, Lang91} for introductions to these subjects.

Let $Y$ be a complex analytic variety. By an \emph{entire curve in $Y$}, we mean a non-constant holomorphic map $\phi:\CC\to Y$. 
The \emph{analytic special set} of $Y$, denoted by $\Sp_{\rm an}(Y)$, is the Zariski closure of the union of the images of all entire curves in $Y$. 
We say that $Y$ is \emph{Brody hyperbolic} (or simply \emph{hyperbolic}) if $\Sp_{\rm an}(Y)=\emptyset$. 

A famous theorem of Brody asserts that a compact complex analytic variety is Brody hyperbolic if and only if it is Kobayashi hyperbolic.
Throughout this paper, ``hyperbolic'' always means ``Brody hyperbolic''.

In the algebraic situation, the counterparts are as follows.
Let $X$ be a projective variety over a field $K$.
The \emph{algebraic special set} of $X$, denoted by $\Sp_{\rm alg}(X)$, is the Zariski closure in $X$ of the union of the images of all non-constant rational maps from abelian varieties over $\overline K$ to $X_{\overline K}$. We say that $X$ is \emph{algebraically hyperbolic} if $\Sp_{\rm alg}(X)=\emptyset$\footnote{It is called ``groupless" in  \cite{JK20, Javanpeykar2021}.}. We note that the terminology ``algebraically hyperbolic" was used by Demailly \cite{Demailly1997} for a different meaning.

Let $X$ be a projective variety over $\CC$. The above two sets of definitions are related by a series of famous conjectures. In fact, Lang conjectures that $\Sp_{\rm alg}(X)=\Sp_{\rm an}(X)$. 
As a special case, Lang's conjecture implies that $X$ is hyperbolic if and only if it is algebraically hyperbolic. 
On the other hand, the Green--Griffiths--Lang conjecture asserts that $X$ is of general type if and only $\Sp_{\rm an}(X)\neq X$. 
We refer to 
\cite[\S VIII.1]{Lang91} for more details on these far-reaching conjectures.

Finally, we formulate the geometric Bombieri--Lang conjecture as follows. 

\begin{conj}[Geometric Bombieri--Lang] \label{main conj}
Let $K$ be a finitely generated field over a field $k$ of characteristic 0 such that $k$ is algebraically closed in $K$. 
Let $X$ be a projective variety over $K$. 
Let $Z$ be the Zariski closure of 
$(X\setminus \Sp_{\rm alg}(X))(K)$ in $X$. 
Then there is a finite set $\{Z_1,\dots, Z_r\}$ of distinct closed subvarieties of $Z$
containing all irreducible components of $Z$ and satisfying the following conditions: 
\begin{enumerate}[(1)]
\item
For each $i=1,\dots, r$, there is a birational $K$-map 
$\rho_i:T_{i,K}\dashrightarrow Z_i$, 
where $T_{i,K}=T_i\times_k {K}$ is the base change for a projective variety $T_i$ over $k$.
\item 
Denote by $U_i$ the maximal open $K$-subvariety of $T_{i,K}$ such that $\rho_i$ extends to a $K$-morphism $\rho_i^\circ:U_i\to Z_i$.
Then the set $(X\setminus \Sp_{\rm alg}(X))(K)$ is contained in the union over $i=1,\dots, r$ of the images of the composition
$$T_i(k) \cap U_i(K) \stackrel{\rho_i^\circ}{\lra}  Z_i(K)\lra Z(K)
\lra X(K).$$
Here the intersection $T_i(k) \cap U_i(K)$ is taken in 
$(T_{i,K})(K)$ via the canonical injection $T_i(k) \to (T_{i,K})(K)$.
\end{enumerate}
Moreover, if $X_{\ol K}$ does not contain any (possible singular) rational curve, then we can take every birational $K$-map 
$\rho_i:T_{i,K}\dashrightarrow Z_i$ 
to be a $K$-morphism under which $T_{i,K}$ is $Z_i$-isomorphic to the normalization $Z_i'$ of $Z_i$. 
\end{conj}

In the conjecture, the last statement applies to the hyperbolic case and the case that $X$ has a finite morphism to an abelian variety. In fact, in these cases, $X_{\ol K}$ does not contain any rational curve. These are the major cases treated in this paper. 
On the other hand, for general $X$, the rational map $\rho_i$ can fail to be a morphism (cf. Example \ref{blow up}). 

In the conjecture, the set $\{Z_1,\dots, Z_r\}$ contains all irreducible components of $Z$, but it may need to contain more elements to cover all the algebraic points as in (2) even in the algebraically hyperbolic case.
This is justified by Example \ref{exam isotrivial}. 

If $X$ is algebraically hyperbolic, with Lang's conjectural equivalence between algebraic hyperbolicity and analytic hyperbolicity, the above conjecture is a precise form to realize Lang's philosophy we quoted at the beginning.
In \S \ref{sec more conjectures}, we formulated various (weaker) versions of Conjecture \ref{main conj}. 
For example, if $k$ is algebraically closed, 
Conjecture \ref{enhanced special set} asserts that the union $\Sp_{\rm alg+const}(X)$ of $\Sp_{\rm alg}(X)$ with all positive-dimensional birationally constant closed subvarieties of $X$ 
is Zariski closed in $X$;
then Conjecture \ref{main conj variant} asserts that 
 $(X\setminus \Sp_{\rm alg+const}(X))(K)$ is finite.

% Moreover, if we only take the positive-dimensional elements of  $\{Z_1,\dots, Z_r\}$, then the conclusion in (2) becomes that all but finitely many $K$-points of $X\setminus \Sp_{\rm alg}(X)$ come from the $k$-points. This appears to be closer to Lang's philosophy and the Bombieri--Lang conjecture over number fields.

% In the last statement, if $Z_i$ is of general type, which is a consequence of the Green--Griffiths--Lang conjecture, then the constant structure $(T_i,\rho)$ for $Z_i$ in (1) is unique up to isomorphism (cf. Theorem \ref{unique constant structure}). 

% A special case of the conjecture happens when $X=Z$; i.e. $(X\setminus \Sp_{\rm alg}(X))(K)$ is Zariski dense in $X$. Assuming the Green--Griffiths--Lang conjecture, this special case actually implies the general case by successively applying it to the irreducible components of the Zariski closures of the rational points.

As recalled above, Conjecture \ref{main conj} (or its suitable variant) is proved 
for curves by \cite{Manin63,Grauert65}, for
subvarieties of abelian varieties by  \cite{Raynaud83, Buium92, Hrushovski96},   
and for smooth projective varieties with ample cotangent bundles by
\cite{Nog82,MD84}.
Moreover, the conjecture for base changes of hyperbolic varieties from $k=\CC$ to a finitely generated extension $K$ is proved by Noguchi \cite[Cor. 4.2]{Nog92}, as an example of his 
hyperbolically embedded case.
We also refer to \cite{Bogomolov1977} and \cite{Javanpeykar2021a} for more results related to constant varieties (i.e., base changes of varieties from $k$ to $K$).

\subsection{Main results}

To state our main results, we will start with the case of hyperbolic varieties finite over abelian varieties.

\subsubsection*{Hyperbolic covers of abelian varieties}

Consider the situation that $X$ has a finite morphism to an abelian variety. In this case, both Lang's conjecture ($\Sp_{\rm alg}(X)=\Sp_{\rm an}(X)$) 
and the Green--Griffiths--Lang conjecture are confirmed, by combining results of Ueno, 
Kawamata and Yamanoi. We refer to Theorem \ref{special set} and Corollary \ref{special set2} for more details. 

The following theorem asserts that the geometric Bombieri--Lang conjecture holds for algebraically hyperbolic varieties with finite morphisms to abelian varieties. 

\begin{thm}[Theorem \ref{hyperbolic main copy}] \label{hyperbolic main}
Let $K$ be a finitely generated field over a field $k$ of characteristic 0 such that $k$ is algebraically closed in $K$. 
Let $X$ be a projective variety over $K$ with a finite morphism $f:X\to A$ for an abelian variety $A$ over $K$. 
Assume that $X$ is algebraically hyperbolic in that $\Sp_{\rm alg}(X)=\emptyset$. 
Then Conjecture \ref{main conj} holds for $X/K/k$.
\end{thm}

Note that we do not assume that $X$ is smooth over $K$ or that $f$ is surjective in the theorem. In particular, $X$ is allowed to be a closed subvariety of $A$, and in this case versions of the geometric Bombieri--Lang conjecture were proved by \cite{Raynaud83, Buium92}. 

In a forthcoming paper, we will prove some case of the theorem without assuming that
$X$ is algebraically hyperbolic.

\begin{remark}
In the case of number fields, if $f:X\to A$ is finite and surjective (but not an isomorphism), 
very little is known in the direction of the Bombieri--Lang conjecture.
See \cite{Corvaja2022} for some results on the sparsity of rational points in this situation.
\end{remark}

\subsubsection*{Partial heights}

Our proofs of the above theorem requires our new notion of partial heights. Let us sketch this  notion and introduce the non-degeneracy conjecture. 

Let $K=\CC(B)$ be the function field of a smooth projective curve $B$ over $\CC$. 
Let $X$ be a projective variety over $K$. 
Let $\pi:\CX\to  B$ be an integral model of $X$ over $B$, which is assumed to be integral, projective and flat. 
For convenience, we identify $X(K)=\CX(B)$, so a point $x\in X(K)$ corresponds to a section $x:B\to \CX$ by abuse of notations.

Let $L$ be a line bundle on $X$, and $\CL$ be a line bundle on $\CX$ extending $L$.
Then we have a height function  
$$h_\CL: \CX(B)\lra \RR, \quad
x\longmapsto  \deg(x^*\CL).$$
The function $h_\CL$ is a Weil height function for $L$. 

Let $\omega_\CL$ be a closed $(1,1)$-form on $\CX$ representing the cohomology class $c_1(\CL)$.
Then we simply have 
$$h_\CL(x)= \int_{B} x^*\omega_\CL,\quad x\in \CX(B).$$
The integration is on the whole section $B$, but partial heights are defined by taking a smaller domain of  integration.

To define a partial height, let $D\subset B$ be a measurable subset, and let $\omega$ be a real $(1,1)$-form on an open neighborhood of the closure of $\pi^{-1}(D)$ in $\CX$. 
Here the measurability makes sense by taking a coordinate chart of $B$, and  ``real'' means the complex conjugate $\overline{\omega}=\omega$.
Define the \emph{partial height function} $h_{(D,\omega)}: X(K)\to \RR$ by
$$
h_{(D,\omega)}(x):=  \int_{D} x^* \omega, \quad x\in \CX(B).
$$

If $D=B$ and $\omega$ represents $c_1(\CL)$, then we recover the original height function $h_\CL$. 
For application of the partial height, we usually take $D$ to be a disc in $B$.

Note that partial heights seem ``smaller'' than the usual Weil height of an ample line bundle. However, the following non-degeneracy conjecture (cf. Conjecture \ref{conj non-deg}) asserts that they are actually comparable.

\begin{conj}[non-degeneracy conjecture] 
Assume that $X$ does not contain any (possibly singular) rational curve.
Let $h_{L}:X(K)\to \RR$ be a Weil height function associated to an ample line bundle 
$L$ on $X$.
Let $h_{(D,\omega)}:X(K)\to \RR$ be a partial height function associated to a pair $(D,\omega)$ on $\CX$, where $D$ is an open disc in $B$. 
Assume that $\omega$ is strictly positive on an open neighborhood of the closure of $\pi^{-1}(D)$ in $\CX$.
Let $\{x_n\}_{n\geq1}$ be a sequence in $X(K)$. 
Assume that $h_{L}(x_n)$ converges to infinity. 
Then $h_{(D,\omega)}(x_n)$ converges to infinity. 
\end{conj}

We will prove that the non-degeneracy conjecture holds for varieties finite over abelian varieties (cf. Theorem \ref{confirm conj non-deg}).

\subsubsection*{The main theorem} 

The following theorem asserts that the non-degeneracy conjecture essentially implies the hyperbolic case of the geometric Bombieri--Lang conjecture.

\begin{thm}[Theorem \ref{partial to BL0 copy}]  \label{partial to BL0}
Let $K=\CC(B)$ be the function field of a smooth projective curve $B$ over $\CC$. 
Let $X$ be a projective variety over $K$. 
Let $\pi:\CX\to B$ be an integral model of $X$ over $B$.
Let $h_{(D,\omega)}:X(K)\to \RR$ be a partial height function associated to a pair 
$(D,\omega)$ on $\CX$, where $D$ is an open subset of $B$ whose closure $\ol D$ is contained in an open disc $D'$ of $B$. 
Assume that the fiber $\CX_b$ is Brody hyperbolic for every $b\in \ol D$.
Assume that the non-degeneracy conjecture holds for $h_{(D,\omega)}$. 
Then Conjecture \ref{main conj} holds for $X/K/\CC$. 
\end{thm}

The hyperbolicity of $\CX_b$ for every $b\in \ol D$ implies the algebraic hyperbolicity of $\CX_b$ for every $b\in \ol D$, and thus implies 
the algebraic hyperbolicity of $X$. 

On the other hand, by a consequence of the Brody lemma (cf. \cite[III, Prop. 3.1]{Lang87}), being hyperbolic is an open condition under the Euclidean topology.
Thus the assumptions of the theorem is satisfied if there is a closed point $b\in B$
such that $\CX_b$ is hyperbolic and that the non-degeneracy conjecture holds for $h_{(D,\omega)}$ for all sufficiently small open neighborhoods of $b$ in $B$. 

Since the non-degeneracy conjecture holds for varieties finite over abelian varieties (cf. Theorem \ref{confirm conj non-deg}), Theorem \ref{hyperbolic main} is essentially a consequence of Theorem \ref{partial to BL0} (or its proof).

\subsection{Ideas of proofs} 

Here we describe our ideas to prove Theorem \ref{hyperbolic main} and Theorem \ref{partial to BL0}. 
Our proofs apply the classical Brody lemma crucially to construct entire curves from sequences of holomorphic sections on the integral model $\CX\to B$. In the following, we will review the Brody lemma, and then outline our idea.
We will give a precise statement of the Brody lemma so that we can refer back from the main body of this article. 

\subsubsection*{The Brody Lemma} 

For compact complex analytic spaces, the equivalence between the Kobayashi hyperbolicity and the Brody hyperbolicity can be established via the Brody lemma. 
In fact, the Brody lemma is a general procedure to construct an entire curve on an analytic variety by a limit process from a suitable sequence of holomorphic maps from discs to the analytic variety.

Denote the standard disc 
$\DD_{r}=\{z\in \CC: |z|<r\}$, and write $\DD=\DD_1$ for simplicity.
Denote by $\ds v_{\rm st}=\frac{\d}{\d z}$ the tangent vector of $\DD_r$ at $0$ under the standard coordinate $z$.

Let $Y$ be a complex analytic variety. Let $\{\phi_n: U_n\to Y\}_{n\geq 1}$ be a sequence of holomorphic maps from Riemann surfaces $U_n$ to $Y$.
By \emph{a re-parametrization of $\{\phi_n: U_n\to Y\}_{n\geq 1}$}, we mean a sequence 
$$\{\phi_n\circ p_n: \DD_{r_n}\to Y\}_{n\geq 1},$$ 
where 
$\{r_n\}_{n\geq 1}$ is a sequence of positive real numbers,  and $\{p_n: \DD_{r_n} \to U_{n}\}_{n\geq 1}$ is a sequence of holomorphic maps.

Let $Y$ be a complex analytic variety endowed with a metric. Let $\{\phi_n: \DD_{r_n}\to Y\}_{n\geq 1}$ be a sequence of holomorphic maps from discs to $Y$.
We say that an entire curve $\phi:\CC\to Y$ is a \emph{limit of $\{\phi_n\}_n$} if there is a subsequence 
$\{n_i\}_{i\geq 1}$ of $\{n\}_{n\geq 1}$, such that 
$\{r_{n_i}\}_{i\geq 1}$ converges to infinity, and such that $\{\phi_{n_i}: \DD_{r_{n_i}}\to Y\}_{i\geq 1}$ converges to $\phi:\CC\to Y$ uniformly on every compact subset of $\CC$.
The last convergence means that for every compact subset $\Omega$ of $\CC$, the restriction $\{\phi_{n_i}|_{\Omega}: \Omega\to Y\}_{i\geq i_0}$ converges to $\phi|_\Omega:\Omega\to Y$ uniformly on $\Omega$, where $i_0$ is a positive integer such that $\Omega \subseteq \DD_{r_{n_i}}$ for every $i\geq i_0$.

Then we have the following Brody lemma. 

\begin{thm} [Brody lemma]   \label{brody}
Let $Y$ be a compact complex analytic variety endowed with a K\"ahler metric. 
Let $\{\phi_n: \DD\to Y\}_{n\geq1}$ be a sequence of holomorphic maps. 
Assume that $\|(\d \phi_n)(v_{\rm st})\|$ converges to infinity. 
Then there exists an entire curve $\phi:\CC\to Y$ which is a limit of a re-parametrization of $\{\phi_n: \DD\to Y\}_{n\geq1}$. 
\end{thm}

Here $\d \phi_n: T_0\DD\to T_{\phi_n(0)}Y$ is the induced map on the tangent spaces. The norm on $T_{\phi_n(0)}Y$ is the one induced by the K\"ahler metric. We refer to \S\ref{sec preliminary} for these notations on singular complex analytic varieties.
We refer to \cite{Duval2017} for a proof the Brody lemma, which also holds in the singular case.

\subsubsection*{Theorem \ref{partial to BL0}} 

Let us first describe the idea to prove Theorem \ref{partial to BL0}. 
By induction, one can reduce it to the case that $X(K)$ is Zariski dense in $X$. 
By a result of Noguchi, it suffices to prove that $X$ is dominated by the base change $V_K$ of a projective variety $V$ from $\CC$ to $K$.
Let $h_L:X(K)\to \RR$ be a Weil height function associated to an ample line bundle $L$ on $X$. 
If $h_L$ is bounded on $X(K)$, then the sections in $\CX(B)$ form a bounded family in the setting of Hilbert schemes of subschemes of $\CX$.
By many geometric operations, the Hilbert scheme eventually gives the variety $V$.

The hard part is to prove that the height function $h_L$ is bounded on $X(K)$. 
Assume the contrary that $h_L$ is unbounded on $X(K)$. 
By the non-degeneracy conjecture, $h_{(D,\omega)}$ is also unbounded on $X(K)$. 
Take an infinite sequence $\{x_n\}$ in $\CX(B)$ with $h_{(D,\omega)}\to \infty$. 
Then we get a sequence  $\phi_n=x_n|_D: D\to \CX$ of maps from the disc to $\CX$.
We plan to apply the Brody lemma to the sequence, by suitably identifying $D$ with the standard unit disc. 
The key is that the property $h_{(D,\omega)}(x_n)\to \infty$ guarantees the condition
 $\|(\d \phi_n)(v_{\rm st})\|\to \infty$, the key condition in the Brody lemma. 
 This can be done by a quick local computation. 
 Therefore, we obtain an entire curve $\phi:\CC\to \CX$, which is a limit of a re-parametrization of $\{\phi_n: D\to \CX\}_{n\geq1}$. 
The entire curve is necessarily contained in a fiber of $\CX_{\ol D}\to \ol D$. 
This contradicts the hyperbolicity assumption of the fibers, and  proves Theorem \ref{partial to BL0}. 

%\newpage
%
%During a discussion with Boucksom, he told us that a similar idea of making horizontal sections vertical  is also used  in Campana's work \cite{Campana1992}.
%More precisely, let $M$ be a hyper-K\"ahler manifold.
%By Yau's theorem \cite{Yau1978}, $M$ admits a hyper-K\"ahler metric and an associated family of complex structures $\tau$ parametrized by $\P^1(\C)$, extending the original one.
%The twistor space $ Z\to \P^1(\C)$ is the complex manifold  obtained by gluing the induced complex structures into a holomorphic family over $\P^1(\C).$
%Campana showed that some fiber $M_{\tau}$ of the  twistor space is not hyperbolic. In his proof, the required entire curve arises as a degenerate deformation of horizontal section of $ Z\to \P^1(\C).$  The non-hyperbolicity of all hyper-K\"ahler manifolds was proved by Verbitsky \cite{Verbitsky2015}. Campana's result plays a key role in Verbitsky's proof.

During a discussion, Boucksom told us that a similar idea of constructing vertical entire curves from horizontal sections was previously used  in Campana's work \cite{Campana1992}.
More precisely, for a hyper-K\"ahler manifold $M$, 
Campana constructed an entire curve on some deformation of $M$, and concluded that this deformation is not hyperbolic.  The non-hyperbolicity of all hyper-K\"ahler manifolds was proved by Verbitsky \cite{Verbitsky2015}, where Campana's result plays a key role.

\subsubsection*{Theorem \ref{hyperbolic main}} 

For Theorem \ref{hyperbolic main}, 
by the Lefschetz principle and some geometric arguments, we can reduce it to the essential case that $k=\CC$ and $K=\CC(B)$ for a smooth projective curve $B$ over $\CC$. 
The key is that in this case the non-degeneracy conjecture holds for $X$. 
Then Theorem \ref{hyperbolic main} is essentially a consequence of Theorem \ref{partial to BL0} (or variants of its proofs). 

Now we describe our proof of the non-degeneracy conjecture for $X$ in Theorem \ref{confirm conj non-deg}.  
By pull-back via the finite morphism $f:X\to A$, it
 suffices to prove the non-degeneracy conjecture for the abelian variety $A$ over $K=\CC(B)$.
The new ingredient here is the canonical partial height. 

In fact, let $U$ be an open subvariety of $B$ such that $A$ has a smooth integral model $\CA_U\to U$ over $U$. 
In the pair $(D, \omega)$ defining the partial height $h_{(D, \omega)}$, take $D$ to be an open disc in $U$, and take $\omega$ to be the Betti form on $\CA_U$ associated to a symmetric and ample line bundle $L$ on $A$. 
The Betti form is a semipositive $(1,1)$-form on $\CA_U$, which represents the cohomology classes of $L$ on fibers of $\CA_U\to U$ and satisfies the nice dynamical property $[m]^*\omega=m^2\omega$ for $m\in \ZZ$. 
The partial height $h_{(D, \omega)}:A(K)\to \RR$ behaves like the N\'eron--Tate height. 
In particular, its induces a positive definite quadratic form on the vector space
$(A(K)/A^{(K/\CC)}(\CC))\otimes_\ZZ\RR$, where $A^{(K/\CC)}$ is the $K/\CC$-trace of $A$.
The vector space is finite-dimensional by the Lang--N\'eron theorem. 
As a consequence, any two positive definite quadratic forms on it can bound each other (up to positive multiples). Take these two quadratic forms to be the partial canonical height $h_{(D,\omega)}$ and the canonical height $\hat h_L$ associated to $L$.

\subsection{Notations and terminology}

For any abelian group $M$ and any ring $R$ containing $\ZZ$, denote
$M_R=M\otimes_\ZZ R$. This apply particularly to $R=\QQ,\RR,\CC$.

By a \emph{variety}, we mean an integral scheme, separated of finite type over the base field. A \emph{curve} is a 1-dimensional variety. 

By a \emph{function field of one variable} over a field $k$, we mean a finitely generated field $K$ over $k$ of transcendence degree 1 such that $k$ is algebraically closed in $K$. 
We usual denote by $B$ a smooth quasi-projective curve over $k$ with function field $K$.
For a projective variety $X$ over $K$, an \emph{integral model of $X$ over $B$}
is a quasi-projective variety $\CX$ over $k$ together with a projective and flat morphism $\CX\to B$ whose generic fiber is isomorphic to $X$.

Let $K$ be a finitely generated field over a field $k$, and assume that $k$ is algebraically closed in $K$.
Let $A$ be an abelian variety over $K$. 
Denote by $A^{(K/k)}$ Chow's $K/k$-trace of $A$ over $k$.
Denote
$$V(A,K):=(A(K)/A^{(K/k)}(k))\otimes_{\ZZ}\RR,$$ 
which is a finite dimensional $\R$-vector space.
In our setting, $k$ is usually $\CC$.

By a \emph{point} of a variety over $\CC$, we mean a closed point. 
By the \emph{generic point} of an integral variety, we mean the generic point of the scheme.

All complex analytic varieties are assumed to be reduced and irreducible.
For a complex analytic variety $X$ with a point $x\in X$, denote by $T_xX$ the complex analytic tangent space of $X$ at $x$ defined by holomorphic derivations. 
For a holomorphic map $f:X\to Y$ with $f(x)=y$, denote by $\d f: T_xX\to T_yY$ the induced map between the tangent spaces.

% For a complex abelian variety $A$, the complex Lie algebra $\Lie(A)$ is defined to be the group of translation-invariant holomorphic derivations on $\CO_A$. By restriction, we have canonical isomorphisms $\Lie(A)\simeq T_xA$ for any $x\in A$.

Let $X$ be a complex analytic space, and $S$ be a subset of $X$. 
Denote by $\ol S$ the closure of $S$ in $X$ under the Euclidean topology. 
Assume that the smooth locus $X^{\mathrm{sm}}$ is covered by countably many open balls $\{U_\alpha\}_{\alpha\in I}$. 
We say that $S$ is \emph{measurable} if $S\cap D_\alpha$ is measurable under the Lebesgue measure of the ball $D_\alpha$ for all 
$\alpha$.  We say that $S$ has measure zero if $S\cap D_\alpha$ has measurable zero for all $\alpha$.

For any positive real number $r$, denote by
$$
\DD_r:=\{z\in \CC: |z|< r\},\qquad
\overline\DD_r:=\{z\in \CC: |z|\leq r\}
$$
the discs of radius $r$. 
Write $\DD=\DD_1$ and $\overline\DD=\overline\DD_1$.
Denote by $\ds v_{\rm st}=\frac{\d}{\d z}$ the tangent vector of $\DD_r$ at $0$ under the standard coordinate $z$.

Let $(Y, d)$ be a metric space. Let $\{r_n\}_{n\geq 1}$ be a sequence of positive real numbers convergent to infinity. Let $\phi_n:\D_{r_n}\to Y$ be a sequence of continuous maps.
We say that $\{\phi_n\}_n$ \emph{converges} to a map $\phi: \C\to Y$ if it converges on every compact subset $\Omega$ of $\C$. Since $\Omega\subseteq \D_{r_n}$ for $n$ sufficiently large, the above definition makes sense. If such $\phi$ exists, it is unique and continuous. Moreover, if $Y$ is further a complex analytic variety and $\phi_n$ are holomorphic, then $\phi$ is holomorphic. 
%We say that $\phi_n$ \emph{uniformly} converges  to $\phi: \C\to Y$ if for every $\epsilon>0$, there is $N\geq 1$ such that for every $n\geq N$ and $z\in \D_{r_n}$, $d(\phi(z),\phi_n(z))\leq \epsilon.$ If $\phi_n$ uniformly converges to $\phi$, then it converges to $\phi.$ 
%One has the following Cauchy convergence criterion:
%$\phi_n$ uniformly converges to some $\phi$ if and only if for every $\epsilon>0$ there is $N\geq 1$ such that for every $m,n\geq N$ and $z\in \D_{\min\{r_n.r_m\}}$ 
%$d(\phi_m(z),\phi_n(z))\leq \epsilon.$

\medskip

\noindent\textbf{Acknowledgment.}
The authors are grateful to Gang Tian, Songyan Xie, and Chenyang Xu for many important communications, and to Shou-Wu Zhang for motivating them to formulate other versions of the geometric Bombieri--Lang conjecture.
The authors would like to thank S\'ebastien Boucksom for telling us Campana's work \cite{Campana1992}. 
The authors would like to thank Ariyan Javanpeykar and Junjiro Noguchi for their helpful comments on the first version of this work.

The authors would like to thank the support of the China-Russia Mathematics Center during the preparation of this paper. 
The first author is supported by a grant from the National Science Foundation of China (grant NO.12271007).
The second author is supported by a grant from the National Science Foundation of China (grant NO. 12250004) and the Xplorer Prize from
the New Cornerstone Science Foundation.

\section{Partial heights on projective varieties} \label{sec partial0}

In this section, we introduce partial heights on projective varieties over function fields of one variable over $\CC$, and introduce its basic properties.
Our key ingredient is a non-degeneracy conjecture of partial heights (Conjecture \ref{conj non-deg}), whose truth implies a weaker version of the geometric Bombieri--Lang conjecture in the hyperbolic case (cf. Theorem \ref{partial to BL0}).
Then in \S \ref{sec more conjectures}, we formulate more versions of the geometric Bombieri--Lang conjecture (cf. Conjecture \ref{main conj}).

\subsection{Partial heights} \label{sec partial00}

The goal of this subsection is to introduce partial heights on projective varieties over function fields of one variable over $\CC$.

\subsubsection*{Differential forms and metrics} 

We review the notion of differential forms on possibly singular complex analytic varieties. Our exposition essentially follows \cite[\S1.1]{King71}.

Recall that a \emph{$(p,q)$-form} $\alpha$ on a (reduced and irreducible) complex analytic variety $X$ is a $(p,q)$-form on the smooth locus $X^{\reg}$ of $X$ such that for every point $x\in X$,  there exists an open neighborhood $U$ of $x$ and a closed embedding $U\hookrightarrow V$ into a complex manifold $V$ such that there exists a $(p,q)$-form $\alpha_V$ on $V$ such that $\alpha_V|_{U\cap X^{\reg}}=\alpha|_{U\cap X^{\reg}}.$ 
Note that the condition is automatic for $x\in X^\reg$ by taking $U=V=X^\reg$, but it is essential at singular points of $X$. 

A $(1,1)$-form $\alpha$ on $X$ is called \emph{K\"ahler} (resp. \emph{positive, semipositive, real}) if 
%it is a K\"ahler form on $X^{\reg}$ and 
for every point $x\in X$, the $(1,1)$-form $\alpha_V$ above can be chosen to be a K\"ahler (resp. positive definite, positive semi-definite, real) form on a suitable $V.$
To emphasize, we may also write ``positive'' as ``\emph{strictly positive}''.

Let $\alpha$ be a $(1,1)$-form $\alpha$ on $X$. 
The for every $x\in X$, $\alpha$ defines a 
\emph{hermitian metric} $\|\cdot\|_\alpha$
on the tangent space $T_xX$. 
This is automatic if $X$ is smooth at $x$, and extended to the singular case by definition. If $\alpha$ is only semipositive, then we still get a semi-norm $\|\cdot\|_\alpha$ on the tangent space $T_xX$. 

The following result is well-known, but we provide a proof for the lack of a direct reference.

\begin{lem} \label{positivity}
Let $X$ be a complex analytic variety, and $W$ be an open subset of $X$ which is contained in a compact subset of $X$. 
Let $\alpha$ be a positive $(1,1)$-form on $X$ and 
$\beta$ be a real $(1,1)$-form on $X$.
Then there is a real constant $c>0$ such that $\alpha|_W-c \beta|_W$ is a positive $(1,1)$-form on $W$.
\end{lem}
\begin{proof}
By compactness, it suffices to prove that for every $x\in X$, there is $c>0$ such that $\alpha-c \beta$ is positive on an open neighborhood of $x$. 
Let $(U, V, \alpha_V)$ be a triple as above representing $\alpha$ at $x$, and 
let $(U', V', \beta_{V'})$ be a triple as above representing $\beta$ at $x$. 
We can assume that $\alpha_V$ is positive on $V$. 
By definition of analytic varieties, replacing $(U, V)$ by their open subsets (containing $x$) if necessary, there are compatible holomorphic maps $U\to U'$ and $V\to V'$. Denote by $\beta_{V}$ the pull-back of $\beta_{V'}$ via $V\to V'$. Then $(U, V, \beta_V)$ is also a triple representing $\beta$ at $x$.
Now the existence of $c$ is easy. 
\end{proof}

\subsubsection*{Partial heights} 

Now we are ready to introduce of partial heights. We first review the usual Weil height. 

Let $K=\CC(B)$ be the function field of a smooth projective curve $B$ over $\CC$. %Note that we have assumed that $B$ is projective for convenience. 
Let $X$ be a projective variety over $K$, and let $L$ be a line bundle on $X$. There is a Weil height function
$$h_L: X(\overline K)\lra \RR$$ 
associated to $L$, which is unique up to bounded functions.

Let $\pi:\CX\to  B$ be an integral model of $X$ over $B$, which is assumed to be integral, projective and flat. 
Let $\CL$ be a line bundle on $\CX$ extending $L$.
Then we have a height function  
$$h_\CL: X(\ol K)\lra \RR, \quad
x\longmapsto \frac{1}{\deg(x)} \deg(\CL|_{\wt x}),$$
where $\wt x$ is the Zariski closure of $x$ in $\CX$.
The function $h_\CL$ is a Weil height function for $L$. 

Let $\omega_\CL$ be a closed $(1,1)$-form on $\CX$ representing the cohomology class $c_1(\CL)$.
Note that if $\CX$ is singular, this condition makes sense by pull-back via a desingularization of $\CX$. 
Then we simply have 
$$h_\CL(x)=\frac{1}{\deg(x)} \int_{\wt x} \omega_\CL,\quad x\in X(\ol K).$$
The integration is on the whole multi-section $\tilde x$, but partial heights are defined by contracting the domain of the integration.

To define a partial height, we need a pair $(D, \omega)$. 
Let $D\subset B$ be a measurable subset. Here the measurability makes sense by taking a coordinate chart of $B$.
Denote by $\ol D$ the closure of $D$ in $B$, which is compact. 
Denote $\CX_D=\pi^{-1}(D)$ and $\CX_{\ol D}=\pi^{-1}(\ol D)$.
Let $\omega$ be a real $(1,1)$-form on an open neighborhood of $\CX_{\ol D}$ in $\CX$. Here ``real'' means the complex conjugate $\overline{\omega}=\omega$.
Define \emph{the partial height of a section $x\in X(\ol K)$ with respect to $(D,\omega)$} to be
$$
h_{(D,\omega)}(x):=\frac{1}{\deg(x)} \int_{\wt x\cap \CX_D} \omega.
$$
This gives a \emph{partial height function} 
$$h_{(D,\omega)}: X(\ol K)\lra \RR.$$
For convenience, we say the pair $(D, \omega)$ is \emph{strictly positive} (resp. \emph{semipositive}) if 
$\omega$ is strictly positive (resp. semipositive) on an open neighborhood of $\CX_{\ol D}$ in $\CX$.

If $D=B$ and $\omega$ represents $c_1(\CL)$, then we recover the original height function. To emphasize, the original height is also called \emph{a full height}. 
For application of the partial height, we usually take $D$ to be a disc in $B$ under a local coordinate.

For our interest, we will only consider the height functions on $X(K)$.  
We make an identification $\CX(B)= X(K)$ via the canonical isomorphism, 
so a point $x\in X(K)$ also represents a section $x:B\to \CX$. 
In this case, we have 
$$
h_\CL(x)=\int_{B} x^*\omega_\CL
$$
and
$$
h_{(D,\omega)}(x)=\int_{D} x^*\omega.
$$

In general, we have the following easy properties.

\begin{lem} \label{comparison partial}
Let $h_{(D,\omega)}$  
be a partial height function on $X(\ol K)$ defined by a pair 
$(D, \omega)$ on $\CX$.
Then the following hold:
\begin{enumerate}[(1)]
\item If $h_{(D, \omega')}$ is a partial height function on $X(\ol K)$ defined by a strictly positive pair $(D, \omega')$ on $\CX$, then there is a constant $c>0$ such that 
$$
-c\, h_{(D, \omega')}\leq  h_{(D,\omega)} \leq c\, h_{(D, \omega')}
$$
on $X(\ol K)$. 

\item If $h_{L}$ is Weil height function on $X(\ol K)$ defined by an ample line bundle 
$L$ on $X$, then there are constants $c_1,c_2>0$ such that 
$$
-c_1 h_{L}-c_2\leq  h_{(D,\omega)} \leq c_1 h_{L}+c_2
$$
on $X(\ol K)$.
\end{enumerate}
\end{lem}
\begin{proof}
For (1), note that Lemma \ref{positivity}
implies that there is a constant $c$ such that 
$-c \omega' \leq   \omega  \leq c\omega'$
on an open neighborhood of $\CX_{\ol D}$. 

To prove (2), we can assume that $h_L=h_\CL$ for an ample line bundle $\CL$ on 
$\CX$ extending $L$, and by (1), we can further assume that $\omega$ is the restriction to $\CX_U$
of a strictly positive $(1,1)$-form on $\CX$ representing $c_1(\CL)$. 
\end{proof}

\subsection{Non-degeneracy}

The goal of this subsection is to present a non-degeneracy conjecture of the partial height, and prove that the conjecture almost settles the geometric Bombieri--Lang conjecture in the hyperbolic case.

\subsubsection*{The non-degeneracy conjecture}

Concerning the growth of the partial heights, we make the following non-degeneracy conjecture. We will see in Theorem \ref{partial to BL0} that the conjecture almost settles the geometric Bombieri--Lang conjecture in the hyperbolic case.

\begin{conj}[non-degeneracy conjecture] \label{conj non-deg}
Let $K=\CC(B)$ be the function field of a smooth projective curve $B$ over $\CC$. 
Let $X$ be a projective variety over $K$. 
Assume that $X$ does not contain any (possible singular) rational curve.
Let $\pi:\CX\to B$ be an integral model of $X$ over $B$.
Let $h_{L}:X(K)\to \RR$ be a Weil height function associated to an ample line bundle 
$L$ on $X$.
Let $h_{(D,\omega)}:X(K)\to \RR$ be a partial height function associated to a strictly positive pair 
$(D,\omega)$ on $\CX$, where $D$ is an open disc in $B$. 
Let $\{x_n\}_{n\geq1}$ be a sequence in $X(K)$. 
Assume that $h_{L}(x_n)$ converges to infinity. 
Then $h_{(D,\omega)}(x_n)$ converges to infinity. 
\end{conj}

By Lemma \ref{comparison partial}(1), the truth of the conjecture does not depend on the choice of $\omega$ in the strictly positive pair $(D,\omega)$. 
We will see from Theorem \ref{confirm conj non-deg} that the conjecture holds for abelian varieties, and thus it also holds for varieties with finite morphisms to abelian varieties.

Note that the conjecture assumes that $X$ does not contain any rational curve. 
This is necessary by the following example, which shows that Conjecture \ref{conj non-deg} fails for $X=\PP^1$.

\begin{example} \label{counter}
Set 
$B=\PP^1_\CC$ and $K=\CC(z)$.  
Set $X=\PP^1_K$ and $\CX=\PP^1_B= B\times_\CC \PP^1_\CC$. 
Denote by $\pi:\CX\to B$ and $f:\CX\to \PP^1_\CC$ the projection morphisms.  
Let $z$ be the standard affine coordinate function on $B$, and $t$ be the standard affine coordinate function on the fibers of $\pi:\CX\to B$.
So $(z,t)$ is an affine coordinate of $\CX$. 
Take $x_n\in X(K)$ to be the morphism 
$$x_n:B\lra \PP^1_B, \quad z\longmapsto (z, a_nz^n).$$
Here $a_n$ will be a positive constant converging to 0.
For the line bundle $\CL=f^*\CO(1)$, we simply have $h_\CL(x_n)=n$ converging to infinity. 
On the other hand, take the K\"ahler form 
$$
\omega= \pi^*\frac{\i \d z\wedge\d \bar z}{(1+|z|^2)^2}+ f^*\frac{\i \d t\wedge\d \bar t}{(1+|t|^2)^2}
$$
to be the sum of the Fubini-Study metrics.
Now it is easy to have 
$$
x_n^*\omega= \frac{\i \d z\wedge\d \bar z}{(1+|z|^2)^2}+ 
\frac{n^2a_n^2 |z|^{2(n-1)}\i\d z\wedge\d \bar z}{(1+a_n^2 |z|^{2n})^2}.
$$
Take $D=\DD_r$ to be the standard disc of radius $r$ in $B$.
A calculation in terms of polar coordinates give
$$
h_{(D,\omega)}(x_n)=\int_D x_n^*\omega
=\frac{2\pi r^2}{1+r^2}
+\frac{2\pi n a_n^2 r^{2n}}{1+a_n^2 r^{2n}}.
$$  
Take $a_n=n^{-n}$. Then the partial height is bounded as $n\to\infty$. 
\end{example}

\subsubsection*{Consequence on the geometric Bombieri--Lang conjecture} 

The following main result of this subsection asserts that assuming Conjecture \ref{conj non-deg}, we can prove a version of the geometric Bombieri--Lang conjecture in the hyperbolic case.

\begin{thm}[Theorem \ref{partial to BL0}]  \label{partial to BL0 copy}
Let $K=\CC(B)$ be the function field of a smooth projective curve $B$ over $\CC$. 
Let $X$ be a projective variety over $K$. 
Let $\pi:\CX\to B$ be an integral model of $X$ over $B$.
Let $h_{(D,\omega)}:X(K)\to \RR$ be a partial height function associated to a pair 
$(D,\omega)$ on $\CX$, where $D$ is an open subset of $B$ whose closure $\ol D$ is contained in an open disc $D'$ of $B$. 
Assume that the fiber $\CX_b$ is Brody hyperbolic for every $b\in \ol D$.
Assume that Conjecture \ref{conj non-deg} holds for $h_{(D,\omega)}$. 

Then Conjecture \ref{main conj} holds for $X/K/\CC$. 
Namely, let $Z$ be the Zariski closure of 
$X(K)$ in $X$. 
Then there is a finite set $\{Z_1,\dots, Z_r\}$ of distinct closed subvarieties of $Z$
containing all irreducible components of $Z$ and satisfying the following conditions: 
\begin{enumerate}[(1)]
\item
For each $i$, the normalization $Z_i'$ of $Z_i$ is constant in the sense that
there is a $K$-isomorphism 
$\rho_i:T_i\times_\CC {K}\to Z_i'$ for a projective variety $T_i$ over $\CC$.
\item 
The set $X(K)$ is the union over $i=1,\dots, r$ of the images of the composition
$$T_i(\CC)\lra (T_i\times_\CC {K})(K)\stackrel{\rho_i}{\lra} Z_i'(K) \lra Z_i(K)\lra Z(K)
\lra X(K).$$
\end{enumerate}
\end{thm}

The theorem is a consequence of its special case $Z=X$. For convenience, we state this special case (in a weaker form) as follows. 

\begin{thm} \label{partial to BL00}
Let $K=\CC(B)$ be the function field of a smooth projective curve $B$ over $\CC$. 
Let $X$ be a projective variety over $K$. 
Let $\pi:\CX\to B$ be an integral model of $X$ over $B$.
Let $h_{(D,\omega)}:X(K)\to \RR$ be a partial height function associated to a pair 
$(D,\omega)$ on $\CX$, where $D$ is an open subset of $B$ whose closure $\ol D$ is contained in an open disc $D'$ of $B$. 
Assume that the fiber $\CX_b$ is Brody hyperbolic for every $b\in \ol D$.
Assume that Conjecture \ref{conj non-deg} holds for $h_{(D,\omega)}$. 

Assume that $X(K)$ is Zariski dense in $X$. 
Then the normalization of $X$ is isomorphic to 
the base change $T_K=T\times_k {K}$ for a projective variety $T$ over $k$, and the complement of the image of the composition
$T(k)\to (T_{K})(K)
\to X(K)$
in $X(K)$ 
is not Zariski dense in $X$. 
\end{thm}

We will first prove Theorem \ref{partial to BL0 copy} by Theorem \ref{partial to BL00}, and then prove Theorem \ref{partial to BL00}.

\subsubsection*{Proof of Theorem \ref{partial to BL0 copy} by Theorem \ref{partial to BL00}} 

Let $(X, Z)$ be as in Theorem \ref{partial to BL0 copy}. 
Let $Z_1, \dots, Z_{r_0}$ be the irreducible components of $Z$. 
By assumption, $Z$ is the Zariski closure of $X(K)$ in $X$, so  
$Z_i(K)$ is Zariski dense in $Z_i$ for $i=1,\dots, r_0$. 
Apply Theorem \ref{partial to BL00} to the variety $Z_i$ for $i=1,\dots, r_0$.
As in part (1) of Theorem \ref{partial to BL0 copy}, we obtain $(T_i,\rho_i)$ for $i=1,\dots, r_0$.
Moreover, the complement of $\Im(T_i(\CC)\to Z_i(K))$ in $Z_i(K)$ is not Zariski dense in $Z_i$. 

To get part (2) of Theorem \ref{partial to BL0 copy}, we might need to enhance the set $\{Z_1, \dots, Z_{r_0}\}$. 
In fact, if $\Sigma:=\cup_{1\leq i\leq r_0} \Im(T_i(\CC)\to Z(K))$ is not equal to 
$Z(K)$, denote by $Z^{(1)}$ the Zariski closure of $Z(K)\setminus \Sigma$ in $X$.
Note that every irreducible component of $Z^{(1)}$ is properly contained in an irreducible component of $Z$. 
Apply Theorem \ref{partial to BL00} to the irreducible components $Z_{r_0+1}, Z_{r_0+2}, \dots, Z_{r_1}$ (with $r_1> r_0$)  of $Z^{(1)}$.
Repeat the process. 
Eventually, we have closed subvarieties $Z_{r_0+1}, Z_{r_0+2}, \dots, Z_{r}$ (with $r\geq r_0$) of $Z$, with a
$K$-isomorphism 
$\rho_i:T_i\times_\CC {K}\to Z_i'$ for a projective variety $T_i$ over $\CC$ for $i=r_0+1, \dots, r$, where $Z_i'$ is the normalization of $Z_i$,
such that 
$$Z(K)=\cup_{1\leq i\leq r} \Im(T_i(\CC)\to Z(K)).$$
This proves Theorem \ref{partial to BL0 copy} (and Theorem \ref{partial to BL0}).

\subsubsection*{Constructing entire curves} 

The key to prove Theorem \ref{partial to BL00} is the following theorem. 
The idea is to use sections $x_n\in \CX(B)$ to construct discs in $\CX$, and the discs gets larger and large in $\CX$
if the 
partial height $h_{(D,\omega)}(x_n)$ converges to infinity, 
so Brody's lemma produces an entire curve on a fiber of $\CX\to B$.

\begin{thm}  \label{partial to BL1}
Let $K=\CC(B)$ be the function field of a smooth projective curve $B$ over $\CC$. 
Let $X$ be a projective variety over $K$. 
Let $\pi:\CX\to B$ be an integral model of $X$ over $B$.
Let $h_{(D,\omega)}:X(K)\to \RR$ be a partial height function associated to a pair 
$(D,\omega)$ on $\CX$, where $D$ is an open subset of $B$ whose closure $\ol D$ is contained in an open disc $D'$ of $B$. 
Assume that the fiber $\CX_b$ is Brody hyperbolic for every $b\in \ol D$.
Then the partial height $h_{(D,\omega)}$ is bounded above on $X(K)$. 
\end{thm}

\begin{proof}
Assume that $h_{(D,\omega)}(x_n)$ converges to infinity for a sequence $\{x_n\}_{n\geq1}$ in $X(K)$. 
We need to find a point $b\in \ol D$ such that the fiber $\CX_b$ is not Brody hyperbolic.

By Lemma \ref{comparison partial}(1), we can assume that $\omega$ is a K\"ahler metric on the whole space $\CX$, which does not violate the assumption that $h_{(D,\omega)}(x_n)\to \infty$.  

Let $\ds v_{\rm st}=\frac{\d}{\d z}$ be the holomorphic vector field on $D'$ by the coordinate function $z$ of $D'$. 
By the definition of the metric,  
$$x_n^*\omega= \|(\d x_n)(v_{\rm st})\|_\omega^2\cdot (i\d z\wedge \d \bar z)$$
on $D'$. 
By assumption,
$$h_{(D,\omega)}(x_n)=\int_{D} x_n^*\omega\to \infty.$$
Hence, there is a sequence $\{b_n\}_{n\geq1}$ in $D$ such that 
$$\|(\d x_n)(v_{\rm st})\|_\omega\to \infty.$$

There is a deceasing sequence $\{r_n\}_{n\geq 1}$ of positive real numbers such that the following holds:
\begin{enumerate}
\item[(1)] for any $b\in \OD$, the open disc $D(b,r_1)$ of center $b$ and radius $r_1$ is contained in $D'$;
\item[(2)] $r_n\to 0;$
\item[(3)] $r_n\|(\d x_n)(v_{\rm st})\|_\omega\to \infty.$
\end{enumerate}
Finally, we are ready to apply the Brody lemma in Theorem \ref{brody} to the sequence 
$\{\phi_n: \DD\to \CX\}_{n\geq 1}$ defined by $z\mapsto x_n(b_n+r_nz)$.
As a result, there is an entire curve $\phi:\CC\to \CX$, as a limit of a re-parametrization of the sequence $\phi_n: \DD\to \CX$. 

As $r_n\to 0$, by construction, the image $\phi(\CC)$ lies in $\CX_\OD\subset \CX_{D'}$.
The map $\pi\circ\phi:\CC\to D'$ is a constant map by Liouville's theorem. 
As a consequence, the image of $\phi:\CC\to \CX$ lies in the fiber $\CX_{b}$ for some $b\in \OD$, and thus induces an entire curve on $\CX_{b}$.
%Inside $D'$, there is $r>0$ such that for any $b\in \OD$, the open disc $D(b,r)$ of center $b$ and radius $r$ is compactly contained in $D'$. 
%Replacing the coordinate function $z$ by a positive multiple if necessary, we can assume $r=1$. 
%Finally, we are ready to apply the Brody lemma in Theorem \ref{brody} to the sequence $\{x_n: \DD(b_n,1)\to \CX\}_{n\geq1}$. 
%As a result, there is an entire curve $\phi:\CC\to \CX$, as a limit of a re-parametrization of the sequence $x_n: \DD(b_n,1)\to \CX$. 
%
%By construction, the image $\phi(\CC)$ lies in $\CA_\OD\subset \CA_{D'}$.
%The map $\pi\circ\phi:\CC\to D'$ is a constant map by Liouville's theorem. 
%As a consequence, the image of $\phi:\CC\to \CX$ lies in the fiber $\CX_{b}$ for some $b\in \OD$, and thus induces an entire curve on $\CX_{b}$.
\end{proof}

\subsubsection*{Constancy 1: dominated by a constant family}

Now we  prove Theorem \ref{partial to BL00}. 
Our proof follows a standard argument using Chow varieties (or Hilbert schemes) with some extra efforts.
For lack of an exact reference, we include a detailed proof here.
%There is a general exposition in \cite[X, Thm. 4.2]{MB85}.

Let $(B/\CC, X/K, \CX/B, (D,\omega))$ be as in Theorem \ref{partial to BL00}.
We first prove that $X/K$ is dominated by a constant variety, and all the rational points come from this process. Namely, 
there is a smooth projective scheme $V$ over $\CC$ (i.e. finite disjoint union of smooth projective varieties over $\CC$), together with a dominant $K$-morphism
$$\varphi:V_K\lra X$$ 
such that composition 
$$V(\CC)\lra V_K(K)\lra X(K)$$ 
is surjective. 
Here we denote $V_K=V\times_\CC K$.

Let $\CL$ be an ample line bundle on $\CX$.
Every point $x\in X(K)$ extends to a section $\wt x\subset \CX$ over $B$.  
Note that $h_{\CL}(x)=\deg_\CL(\wt x)$ corresponds to a Weil height
associated to the ample line bundle $\CL_K$ on $X$.
Applying Theorem \ref{partial to BL1},
we see that the partial height $h_{(D,\omega)}(x)$ is bounded above for $x\in X(K)$.
By assumption, Conjecture \ref{conj non-deg} holds for $h_{(D,\omega)}(x)$,  so the height $h_{\CL}(x)$ is bounded above for $x\in X(K)$.
In other words, 
there is a constant $c>0$ such that $\deg_\CL(\wt x)<c$ for every $x\in X(K)$.

The key is that there is a quasi-projective scheme $S$ over $\CC$ parametrizing the set of sections $Y$ of $\CX\to  B$ with $\deg_\CL(Y)<c$. 
In fact, by the theory of Hilbert schemes or Chow varieties (cf. \cite[Lem. 2.4]{Grothendieck95}), there is a quasi-projective scheme $S'$ over $\CC$ parametrizing 1-dimensional reduced and equi-dimensional subschemes $Y'$ of $\CX$ with $\deg_\CL(Y')<c$. 
From $S'$ to $S$, the extra condition that the composition $Y'\to \CX\to  B$ is an isomorphism is equivalent to that the morphism $Y'\to  B$ is finite and flat of degree 1.
The ``flat'' condition is an open condition on $Y'$, and it implies the  
 ``finite'' condition by a dimension reason. 
Under this condition, the ``degree 1'' condition is an open and closed condition on $Y'$.  Then $S$ is quasi-projective over $\CC$. 

By the moduli, we have a $\CC$-morphism $S\times_\CC B\to  \CX$, which sends the fiber $s\times_\CC B$ for $s\in S(\CC)$ to the section of $\CX\to  B$ corresponding to $s$. 
It follows that the composition 
$$S(\CC)\lra (S\times_\CC B)(B)\lra \CX(B)\lra X(K)$$ 
is surjective. 
By assumption, $X(K)$ is Zariski dense in $X$, so the morphism $S\times_\CC B\to  \CX$ is dominant. 
Taking a base change by $\Spec K\to  B$, we obtain a dominant $K$-morphism
$S\times_\CC K\to X$. 

It is easy to adjust the quasi-projective scheme $S$ to a smooth projective scheme over $\CC$.
In fact, by taking the disjoint union of irreducible components, we can assume that $S$ is a disjoint union of quasi-projective varieties.
Let $V'$ be a disjoint union of projective varieties containing $S$ as a dense and open subscheme.
Let $V$ be a resolution of singularities of $V'$ by Hironaka's theorem.
The morphism $S_K\to X$ gives a dominant rational map $V_K\dashrightarrow X$. We claim that the rational map extends to a morphism $V_K\to X$.

By assumption, the fiber $\CX_b$ of $\CX\to B$ above any $b\in D$ is hyperbolic. 
Then $\CX_b$ does not contain any rational curve.
Consequently, the generic fiber $X$ does not contain any rational curve.
By \cite[Lem. 3.2, Lem. 3.5]{JK20}, which originally follows from \cite[Prop. 6.2]{GLL15}, the rational map $V_K\dashrightarrow X$ extends to a morphism. 

By construction, the composition
$$V(\CC)\lra (V\times_\CC B)(B)\lra \CX(B)\lra X(K)$$ 
is surjective. 
It follows that the morphism $V_K\to X$ is surjective.
%This proves the first statement of Theorem \ref{partial to BL00}. 

\subsubsection*{Constancy 2: constancy in the normal case}

In the above, we have proved that $X$ is dominated by a constant variety $V\times_\CC K$. The following result asserts that $X$ itself is constant if it is normal. 

\begin{thm} \label{constancy}
Let $K=\CC(B)$ be the function field of a smooth projective curve $B$ over $\CC$. 
Let $X$ be a normal projective variety over $K$. 
Let $\CX\to B$ be an integral model of $X$ over $B$. 
Assume that one of the following conditions holds:
\begin{enumerate}[(1)]
\item for some closed point $b\in B$, the fiber $\CX_{b}$ is hyperbolic;
\item $X$ is a variety of general type over $K$.  
\end{enumerate}
Assume that there is a smooth projective variety $V$ over $\CC$ together with a surjective $K$-morphism $\varphi:V\times_\CC {K} \to X$.
Then $X$ is isomorphic to 
the base change $T\times_\CC {K}$ for a projective variety $T$ over $\CC$.
\end{thm}

We will only need case (2), which is essentially proved by \cite{Nog85}, but we add case (1) here since it follows by a similar proof. 
For the sake of readers, we will sketch a proof of the theorem in the next subsection. The main idea of our proof still follows those of 
\cite[p. 37]{Nog85} and \cite[Thm. 6.9.5]{Kob98}.
%Now let us finish the proof of Theorem \ref{partial to BL00}.

Return to the proof of Theorem \ref{partial to BL00}.
In the theorem, we can assume that $X$ is normal by passing to its normalization.  
In the above, we have already proved that
there is a smooth projective scheme $V$ over $\CC$ together with a dominant $K$-morphism
$\varphi:V_K\to X$ 
such that the composition 
$V(\CC)\to V_K(K)\to X(K)$
is surjective. 
By Theorem \ref{constancy}, $X\simeq T_K$ for a projective variety $T$ over $\CC$. 
This gives the first statement of the theorem. 

It remains to prove that via the isomorphism $X\simeq T_{K}$, the complement of 
$\Im(T(\CC)\to X(K))$ in $X(K)$
is not Zariski dense in $X$. 
We need the condition that the map
$V(\CC)\to X(K)$
is surjective.
For the purpose here, we are not able to assume that $V$ is connected. 
Let $V_1, \dots, V_r$ be all the connected components of $V$ such that  
$V_{i,K}\to X$ is surjective.
Then it suffices to prove that 
$\Im(T(\CC)\to X(K))$
contains $\Im(V_i(\CC)\to X(K))$ for every $i$.
It suffices to prove that the $K$-morphism $V_i\times_\CC K \to X\simeq  T_K$ is equal to the base change of a $\CC$-morphism $V_i \to T$.
This is given by the following result.

\begin{thm} \label{unique constant structure}
Let $K/k$ be a field extension of characteristic 0 such that $k$ is algebraically closed in $K$. 
Let $T, T'$ be projective varieties over $k$.
Assume that one of the following conditions:
\begin{enumerate}[(1)]
\item there exists an embedding $k\hookrightarrow \CC$ via which $T\times_k\CC$ is hyperbolic;
\item $T$ is of general type over $k$.
\end{enumerate}
Then any surjective $K$-morphism $T'_K\to T_K$ descends to a unique surjective $k$-morphism $T'\to T$;  any $K$-isomorphism $T'_K\to T_K$ descends to a unique $k$-isomorphism $T'\to T$. 
\end{thm}

\begin{proof}

The first statement implies the second one. In fact, if a $K$-morphism $T'_K\to T_K$ descends to a $k$-morphism $T'\to T$, then $T'_K\to T_K$ is an isomorphism if and only if $T'\to T$ is an isomorphism. 

Now we prove the first statement. 
In case (2), by the Lefschetz principle, we can assume that $k$ is a finitely generated field over $\QQ$. Then we can also fix an embedding $k\hookrightarrow \CC$, as in case (1).

Denote by $\underline{\Hom}_k(T', T)$ the functor sending any noetherian $k$-scheme $S$ to the set of $S$-morphisms $T'_{S}\to T_{S}$. 
By Grothendieck \cite[\S4]{Grothendieck95}, the functor 
 is represented by a scheme locally of finite type over $k$. 
Denote by $\underline{\rm Sur}_k(T', T)$ the sub-functor 
of $\underline{\Hom}_k(T', T)$ corresponding to surjective $S$-morphisms $T'_{S}\to T_{S}$. 
For any $S$-morphism $T'_{S}\to T_{S}$, the dimensions of the fibers of 
$\Im(T'_{S}\to T_{S})$ over $S$ form an upper semi-continuous function on $S$ (cf. \cite[IV-3, Thm. 13.1.5]{EGA4}).
As a consequence,
$\underline{\rm Sur}_k(T', T)$ is represented by a closed subscheme $\CS$ of $\underline{\Hom}_k(T', T)$.

The key is that the set
$\CS(\CC)$ of surjective morphisms $T'_{\CC}\to T_\CC$ is finite. 
This follows from Noguchi's theorem (cf. \cite[Thm. A]{Nog92}) in case (1) and Kobayashi--Ochiai's theorem (cf. \cite[Thm. 7.6.1]{Kob98}) in case (2).

As the scheme $\CS$ is locally of finite type over $k$, we conclude that it is actually finite over $k$.
By assumption, $k$ is algebraically closed in $K$, so $\CS(k)=\CS(K)$. 
It follows that the $K$-morphism $T'_{K}\to T_K$ descends to a $k$-morphism $T'\to T$. 
\end{proof}

As a consequence, the constant structure $(T_i,\rho)$ for $Z_i'$ in in Theorem \ref{partial to BL0 copy}(1) is unique up to isomorphism.

\subsection{Constancy}

The goal of this subsection is to sketch a proof of Theorem \ref{constancy}.
As mentioned above, the theorem is essentially due to Noguchi \cite{Nog85}. 
The main idea of our proof still follows those of 
\cite[p. 37]{Nog85} and \cite[Thm. 6.9.5]{Kob98}.

Let $(B/\CC, X/K, V/\CC, \CX/B)$ be as in the theorem. 
The $K$-morphism $\varphi:V_K\to X$ extends to a surjective $U$-morphism 
$\varphi:V_U\to \CX_U$ for some non-empty open subvariety $U$ of $B$, where we denote $V_U=V\times_\CC U$ and $\CX_U=\CX\times_B U$.
 
The hard part of the proof is to prove that up to shrinking $U$, 
the family $\CX_U\to U$ is analytically locally trivial, i.e., there is an analytic open covering $\{U_i\}_{i\in I}$ of $U$ such that for every $i\in I$, the induced family $\CX_{U_i}\to U_i$ is biholomorphic to the trivial family $T_i\times U_i\to U_i$ for some complex projective variety $T_i$. 
The main idea to use foliation to construction horizontal leaves of $\CX_U\to U$. 
We first treat the case that $X$ is smooth over $K$, and then modify the idea to treat the normal case.

\subsubsection*{Isotriviality:  smooth case}
Assume that $X$ is smooth over $K$.
We can further assume that $\CX_U\to U$ is smooth by shrinking $U$. The goal here is to prove that $\CX_U\to U$ is analytically locally trivial. 
 
Denote by $T\CX_U$ the tangent sheaf of $\CX_U$ over $\CC$ (as a sheaf), denote by
$\CT=\Spec (\Sym^* ((T\CX_U)^\vee))$ the geometric tangent bundle of $\CX_U$ (as a variety).
Then a complex point of $\CT$ is represented by a pair $(t,x)$ with $x\in \CX_U$ a point and $t\in T_x\CX_U$ a tangent vector of $\CX_U$ at $x$. 
The structure morphism $\CT\to \CX_U$ maps $(t,x)$ to $x$. 

Fix a section $t_0\in \Gamma(U,TU)$ which is non-zero at every point of $U$. Note that $t_0$ exists by shrinking $U$ if necessary. 
Consider the map 
$$
\Phi:V\times_\CC U \lra \CT, \quad
(v,u)\longmapsto ((\d \varphi)(t_{0,v,u}),\varphi(v,u)).
$$
Here $\d\varphi:TV_U\to T\CX_U$ is the morphism between tangent sheaves induced by $\varphi:V_U\to \CX_U$, and 
$t_{0,v,u}$ is the tangent vector of $V_U$ at $(v,u)$ given by $(0,t_0)$ via the natural decomposition $T_{(v,u)}(V\times_\CC U)=T_{v}V\oplus T_{u}U$.
 
Note that $\Phi$ is an algebraic morphism over $\CX_U$. 
As $V_U$ is projective over $\CX_U$, the image $\Phi(V_U)$ is also projective over $\CX_U$. As $\CT$ is affine over $\CX_U$, its closed subvariety $\Phi(V_U)$ is also affine over $\CX_U$. 
As a consequence, the morphism $\delta:\Phi(V_U)\to \CX_U$ is finite. 
We are going to define a vector field $\ft\in\Gamma(\CX_U,T\CX_U)$ in terms of $\delta$ .

We will start with some extra notations due to singularity. 
Denote by $\gamma:W\to \Phi(V_U)$ the normalization of $\Phi(V_U)$. 
Then the morphism $\lambda=\delta\circ \gamma:W\to \CX_U$
is still finite. 
Denote by $W^{\rm sing}$ the singular locus of $W$.
Denote $\CX_U^\circ=\CX_U\setminus \lambda(W^{\rm sing})$ and $W^\circ=W\setminus \lambda^{-1}(\lambda(W^{\rm sing}))$.
%Denote $\Phi(V_U)^\circ=\gamma(W^\circ)$. 
Note that $W$ is normal, so 
$W^{\rm sing}$ has at least codimension two in $W$, and thus
$\CX_U\setminus \CX_U^\circ$ has at least codimension two in $\CX_U$. 
Note that $\lambda^\circ:W^\circ\to \CX_U^\circ$ is a finite surjective morphism of regular varieties, so it is flat by the miracle flatness (cf. \cite[Thm. 23.1]{Matsumura}).

By the fiber product, the $\CX_U$-morphism $W\to \CT$ induces a $W$-morphism $W\to \CT\times_{\CX_U}W$. 
This gives a regular section $\fs$ of $\Gamma(W, \lambda^*T\CX_U)$
for the morphism $\lambda:W\to \CX_U$. 
Define $\ft^\circ$ as the image of $\fs$ under the trace map 
$$
\Gamma(W^\circ, \lambda^*T\CX_U) \lra \Gamma(\CX_U^\circ, T\CX_U).
$$
Here we recall the trace map briefly. 
In general, for a finite and flat morphism $f:A\to B$ of noetherian schemes, and for a locally free sheaf $E$ on $B$, there is a canonical trace map 
$f_*f^*E\to E$ of $\CO_B$-modules, whose global section gives a trace map $\Gamma(A,f^*E)\to \Gamma(B,E)$. 
To define them, by the canonical isomorphism $f_*f^*E\to E\otimes f_*\CO_A$, it suffices to introduce a trace map 
$f_*\CO_{A}\to \CO_{B}$.
Note that $f$ is finite and flat, and thus locally free, so we can assume that $f$ is actually free by taking affine open covers. 
Then the trace map comes from the usual trace of ring extensions. Namely, the trace of a  section $\alpha$ of $f_*\CO_A$ is defined as the trace of the multiplication map $\alpha:f_*\CO_{A}\to f_*\CO_{A}$, viewed as a linear endomorphism of a free $\CO_B$-module.

Hence, we have defined a section $\ft^\circ$ of $T\CX_U^\circ$ over $\CX_U^\circ$.
Since $\CX_U$ is regular and $\CX_U\setminus \CX_U^\circ$ has at least codimension two in $\CX_U$, the section $\ft^\circ$ of $T\CX_U^\circ$  extends to a section $\ft$ of $T\CX_U$ uniquely. 
This finishes the definition of $\ft$.

The key property of $\ft$ is that $(\d\pi)(\ft)=\deg(\lambda) t_0$ holds for the morphism 
$\d\pi: T\CX_U\to TU$.
In particular, $\ft$ is non-vanishing at every point of $\CX_U$.
To check the equality, 
it suffices to check it for the morphism 
$(\d\pi)_x: T_x\CX_U\to T_{\pi(x)}U$ at every point $x\in \CX_U$ over which 
$\delta:\Phi(V_U)\to \CX_U$ is \'etale.
For such an $x$, the fiber $\delta^{-1}(x)$ is reduced, the variety $\Phi(V_U)$ is smooth at $\delta^{-1}(x)$, and $W\to \Phi(V_U)$ is isomorphic at $\delta^{-1}(x)$.
By the construction of the trace map, we see that 
the fiber $\ft_x\in T_x \CX_U$ of $\ft$ at $x$ is given by 
$$
\ft_x=\sum_{(t,x)\in \delta^{-1}(x)} t.
$$
Here elements of $\Phi(V_U)\subset \CT$ takes the form $(t,x)$ as mentioned above.
It follows that
$$
(\d\pi)_x(\ft_x)=\sum_{(t,x)\in \delta^{-1}(x)} (\d\pi)_x(t)=\sum_{(t,x)\in \delta^{-1}(x)} t_0=\deg(\lambda) t_0.
$$
This proves the equality.

In the complex analytic setting, consider the foliation associated to the vector field 
$\ft$ on $\CX_U$. The equality $(\d\pi)(\ft)=\deg(\delta) t_0$ implies that every leaf of $\ft$ is \'etale over $U$.
By the foliation, for any point $x\in \CX_U$, there is an open neighborhood $N$ of $x\in  \CX_U$, which gives an open neighborhood $N_1=N\cap \CX_{\pi(x)}$ of $x\in \CX_{\pi(x)}$, and
an open neighborhood $N_2=\pi(N)$ of $\pi(x) \in U$, such that 
there is a biholomorphic map $N\to N_1\times N_2$ satisfying the following conditions:
\begin{enumerate}[(1)]
\item
the composition 
$N\to N_1\times N_2\stackrel{p_2}{\to} N_2$ is equal to $\pi|_N:N\to N_2$;
\item
 the composition 
$N_1\hookrightarrow N\to N_1\times N_2\stackrel{p_1}{\to} N_1$ is the identity map;
\item
the fibers of the composition 
$N\to N_1\times N_2\stackrel{p_1}{\to} N_1$
 are exactly the maximal leaves of $\ft$ in $N$.
\end{enumerate}

For any point $u\in U$, apply the above decomposition to every $x\in \CX_u$. By compactness, there is an open neighborhood $N_2$ of $u \in U$, together with $N=\pi^{-1}(N_2)$ and $N_1=\CX_{u}$, such that 
there is a biholomorphic map $N\to N_1\times N_2$ satisfying the above three conditions.
Therefore, we have finally proved that the fibration $\pi: \CX_U\to U$ is analytically locally trivial.

\subsubsection*{Isotriviality:  normal case}

Now we consider the case that $X$ is normal. 
By shrinking $U$, we can assume that $\CX_U$ is normal, and that all fibers of $\CX_U\to U$ are normal. The goal here is still to prove that $\pi:\CX_U\to U$ is analytically locally trivial. 

The tangent sheaf $T\CX_U$ (defined by derivations) is generally not locally free, so we will avoid using it.
Denote by $\CX_U^\sing$ and $\CX_U^\circ$ the singular locus and the smooth locus of $\CX_U$ over $\CC$.
The tangent sheaf $T\CX_U^\circ$ of $\CX_U^\circ$ is locally free over $\CX_U^\circ$. 

As in the smooth case, fix a section $t_0\in \Gamma(U,TU)$ which is non-zero at every point of $U$ (up to shrinking $U$). 
By modifying the proof of the smooth case above, we can construct a section $\ft$ of $T\CX_U^\circ$ such that $(\d\pi)_x(\ft)=a t_0$ holds for every point $x\in\CX_U^\circ$ and for some positive integer $a$.

Fix an (algebraic) closed immersion $i:\CX_U\to \PP^n_U$ over $U$ for some positive integer $n$. 
%Denote by $i^\circ:\CX_U^\circ\to \PP^n_U$ the induced morphism. 
Denote $\CY=\PP^n_U$ for simplicity.
Then there is an injection $(\d i)^\circ: T\CX_U^\circ \to (T\CY)|_{\CX_U^\circ}$. 
Here $(T\CY)|_{\CX_U^\circ}$ denotes the pull-back of $T\CY$ to $\CX_U^\circ$ as a locally free sheaf.
Then we obtain a section $\ft'=(\d i)^\circ(\ft)$ of $(T\CY)|_{\CX_U^\circ}$ over 
$\CX_U^\circ$.
As $\CX_U$ is normal, the singular locus $\CX_U^\sing$ has codimension at least two in $\CX_U$, 
and the section $\ft'$ extends to a unique section of the locally free sheaf
$(T\CY)|_{\CX_U}$ over $\CX_U$.
We still denote it by $\ft'$.
The original condition $(\d\pi)_x(\ft)=a t_0$ transfers to a similar condition for $\ft'$, so $\ft'$ is nonzero at any point of $\CX_U$. 

Let $x\in \CX_U$ be a closed point. Then there is an open neighborhood (under the Euclidean topology) $M$ of $x$ in $\CY$ such that $\ft'|_{M\cap \CX_U}$ extends to a section $\tilde\ft$ of $(T\CY)|_M=TM$. 
As $\ft'$ is nonzero at any point of $\CX_U$, we can assume that 
$\tilde\ft$ is nonzero at any point of $M$ by shrinking $M$. 

Consider the foliation on $M$ associated to the section $\tilde\ft$ of $TM$. 
Up to shrinking $M$, we can assume that there is a biholomorphic map $M\to M_1\times M_2$ with $M_1=M\cap \CY_{\pi(x)}$ and $M_2=\pi(M)\subset U$, such that the induced projection $p_2:M\to M_2$ is compatible with $\pi$, 
the composition $M_1\hookrightarrow M\stackrel{p_1}{\to} M_1$ is the identity map, and that all fibers of $p_1:M\to M_1$ are maximal leaves of $\tilde\ft$ on $M$. 
Moreover, after shrinking $M$, we may assume that $M_2$ is a disk.

Denote $N=M\cap \CX_U$, $N^\circ=M\cap \CX_U^\circ$ and $N_1=M_1\cap \CX_U$.
Consider the foliation on $N^\circ$ associated to the section $\ft$ of $TN^\circ$. 
By compatibility, for any point $y\in N^\circ$, any local leaf $\CF(y,\ft)$ of $\ft$ on $N^\circ$ through $y$ is also a local leaf of $\tilde\ft$ on $M$ through $y$. 
The intersection $N \cap p_1^{-1}(p_1(y))$ is an analytic subvariety of the disk $p_1^{-1}(p_1(y))\simeq M_2$, and contains the open subset $F(y, t)$ of  $p_1^{-1}(p_1(y))$. So we must have $N \cap p_1^{-1}(p_1(y))=p_1^{-1}(p_1(y))$, and thus $p_1^{-1}(p_1(y)) \subset N$.
By approximation, $p_1^{-1}(p_1(y))$ is contained in $N$ for any $y\in N$. 

Hence, the image $p_1(N)=p_1(N_1)=N_1$, and the restriction of $M\to M_1\times M_2$ induces a bijective map $N\to N_1\times M_2$. 
Note that $N, N_1, M_2$ are normal, so $N\to N_1\times M_2$ is a biholomorphic map.

In the above, the fiber of the projection $N\to N_1$ through any $y\in N^\circ$ is the Zariski closure of the leaf $\CF(y,\ft)$. 
By approximation, we see that that the germ of the map 
$N\to N_1$ at $x$ is independent of the choice of the lifting $\tilde\ft$ of $\ft$ to $M$. 
The same result holds for the isomorphism $N\to N_1\times M_2$. 
Therefore, when varying $x$, we can glue these maps $N\to N_1\times M_2$ together.
As in the smooth case, for any point $u\in U$, we can glue the maps along an open neighborhood of $\CX_u$ in $\CX_U$.  
As a result,  the fibration $\pi: \CX_U\to U$ is analytically locally trivial.

\subsubsection*{Constancy}

Once the fibration $\pi: \CX_U\to U$ is analytically locally trivial, it is easy to prove that it is actually trivial (both analytically and algebraically). 

In case (1) of Theorem \ref{constancy}, by assumption, $\CX$ has a hyperbolic closed fiber over $B$. 
By \cite[III, Prop. 3.1]{Lang87}, being hyperbolic is an open condition under the Euclidean topology.
Then there is a point $u_0\in U$ such that the fiber $\CX_0=\CX_{u_0}$ is hyperbolic. 
Then $\Aut(\CX_0)$ is finite by Noguchi \cite[Thm. A]{Nog92}.

In case (2) of Theorem \ref{constancy}, by assumption,
$X$ is of general type. Then there is a closed point $u_0\in U$ such that the fiber $\CX_0=\CX_{u_0}$ is integral and of general type.
This simple fact can be proved by considering rational maps from $X$ to projective spaces via global sections of pluri-canonical bundles of $X$.
Then $\Aut(\CX_0)$ is finite by 
 Kobayashi--Ochiai's theorem (cf. \cite[Thm. 7.6.1]{Kob98}).

In both cases, the monodromy action gives a homomorphism
$$
\rho:\pi_1(U,u_0)\lra \Aut(\CX_0).
$$
The quotient map $\pi_1(U,u_0)\to \Im(\rho)$ corresponds to a finite unramified covering $(U',u_0')\to (U,u_0)$. 

We still need the surjective $U$-morphism $\varphi:V\times_\CC U \to \CX_U$. 
By base change, we obtain a surjective $U'$-morphism $\varphi':V\times_\CC U' \to \CX_{U'}$. 
By the trivialization $\iota:\CX_{U'}\to \CX_0\times_\CC U'$, we obtain 
a surjective $U'$-morphism $\psi:V\times_\CC U' \to \CX_0\times_\CC U'$. 
By Theorem \ref{unique constant structure},
 $\psi:V\times_\CC U' \to \CX_0\times_\CC U'$ is equal to the base change of a morphism $\psi_0:V \to \CX_0$.
 
 Any point $v\in V(\CC)$ transfers to horizontal sections 
$\tilde v\in (V\times_\CC U)(U)$, $\tilde v'\in (V\times_\CC {U'})(U')$,
and $\psi(\tilde v')\in \CX_{U'}(U')$. 
By compatibility, the image of $\psi(\tilde v')\in \CX_{U'}(U')$ in $\CX_U$ gives a section of $\CX_U(U)$. 
As a consequence, the monodromy action of $\Aut(U'/U)$ on $\CX_0$, depending on the choice of a point $u_0'\in U'$ above $u_0\in U$ and an identification $\CX_{u_0'}\simeq \CX_0$, fixes the point $\psi_0(v)\in \CX_0$. 
Since $\psi_0:V \to \CX_0$ is surjective, 
$\Aut(U'/U)$ actually 
acts trivially on $\CX_0$. 
As a consequence, $\pi: \CX_U\to U$ is a trivial family. 

Therefore, $\CX_U$ is biholomorphic to $\CX_0\times_\CC U$ over $U$. 
This proves that $\pi: \CX_U\to U$ is analytically trivial (or constant). 

We can use a GAGA-type argument to prove that $\pi: \CX_U\to U$ is actually algebraically trivial. 
In fact, fix a point $u_0\in U$ such that the fiber $\CX_0=\CX_{u_0}$ is hyperbolic as above.
The analytically locally trivial fibration $\pi: \CX_U\to U$ is classified by the \v Cech cohomology
$$
H^1(U, \Aut(\CX_0))\simeq \Hom(\pi_1(U,u_0),\Aut(\CX_0)).
$$
As a consequence of Riemann's existence theorem 
(cf. \cite[XII, Thm. 5.1, Cor. 5.2]{SGA1}), the \'etale fundamental group $\pi_1^{\rm et}(U, u_0)$ is canonically isomorphic to the profinite completion of the topological fundamental group $\pi_1(U,u_0)$. 
Then by \cite[Prop. 5.7.19, Prop. 5.7.20]{Fu2011},  the finiteness of $\Aut(\CX_0)$ implies
$$
\Hom(\pi_1(U,u_0),\Aut(\CX_0)) \simeq  \Hom(\pi_1^{\rm et}(U,u_0),\Aut(\CX_0)) 
\simeq H^1_{\rm et}(U, \Aut(\CX_0)).
$$
The last term is the \'etale \v Cech cohomology, which is defined even if $\Aut(\CX_0)$ is non-abelian (cf. \cite[p. 229]{Fu2011}).
This establishes an equivalence between analytically locally trivial families and \'etale locally trivial algebraic families. Moreover, the triviality is compatible in both cases. 
This finishes the proof of Theorem \ref{constancy}.

\subsubsection*{Alternative approach}

A large part of the above argument is to prove that the fibration $\pi: \CX_U\to U$ is analytically locally trivial using the foliation method. In the case that $\pi: \CX_U\to U$ is smooth, if we further assume that
$\CX_b$ is both hyperbolic and of general type for all $b$ in a disc of $U$, then we have an algebraic proof sketched as follows.  

First, for the surjective morphism $V\times U\to \CX_U$, by taking successive generic hyperplane sections of $V$, we can assume that $V\times U\to \CX_U$ is generically finite (and surjective). Here $V$ is assumed to be connected, and it is of general type since $\CX_b$ is of general type. 

Second, by the resolved Iitaka--Severi conjecture (cf. \cite[Thm. 4.1]{GP}), the set $\{\CX_b:b\in D\}$ falls into finitely many birational equivalence classes. Then for some fixed $b_0\in D$, the fiber $\CX_b$ is birational to $\CX_0=\CX_{b_0}$ for uncountably many $b\in  D$. 

Third, the above $\CX_b$ is actually isomorphic to $\CX_0$. This follows from 
\cite[Cor. 6.3.10]{Kob98} or the algebraic version \cite[Lem. 3.2, Lem. 3.5]{JK20} again. 

Fourth, the functor $\CF:=\underline{\mathrm{Isom}}_U(\CX_0\times U, \CX_U)$ over $U$ is representable by a countable union of schemes of finite type over $U$, as a consequence of the theory of Hilbert schemes (cf. \cite[\S4]{Grothendieck95}). 
Note that the third step implies that the fiber $\CF_b$ is non-empty for un-countably many $b\in D$. 
It follows that the generic fiber $\CF_K$ is also non-empty. 
Any point of $\CF_K(\overline K)$ implies that $X$ is a twist of $\CX_0$, and thus $\pi:\CX_U\to U$ is \'etale locally trivial up to shrinking $U$.

\subsection{Geometric Bombieri--Lang conjectures}  \label{sec more conjectures}

Recall that we have formulated the geometric Bombieri--Lang conjecture 
in Conjecture \ref{main conj}. The goal of this subsection is to present some examples related to the conjecture, and more importantly formulate other versions of this conjecture. We will see that Conjecture \ref{main conj} essentially implies the other versions we will present. 
For convenience, we duplicate Conjecture \ref{main conj} as follows.

\begin{conj}[Conjecture \ref{main conj}] \label{main conj copy}
Let $K$ be a finitely generated field over a field $k$ of characteristic 0 such that $k$ is algebraically closed in $K$. 
Let $X$ be a projective variety over $K$. 
Let $Z$ be the Zariski closure of 
$(X\setminus \Sp_{\rm alg}(X))(K)$ in $X$. 
Then there is a finite set $\{Z_1,\dots, Z_r\}$ of distinct closed subvarieties of $Z$
containing all irreducible components of $Z$ and satisfying the following conditions: 
\begin{enumerate}[(1)]
\item
For each $i=1,\dots, r$, there is a birational $K$-map 
$\rho_i:T_{i,K}\dashrightarrow Z_i$, 
where $T_{i,K}=T_i\times_k {K}$ is the base change for a projective variety $T_i$ over $k$.
\item 
Denote by $U_i$ the maximal open $K$-subvariety of $T_{i,K}$ such that $\rho_i$ extends to a $K$-morphism $\rho_i^\circ:U_i\to Z_i$.
Then the set $(X\setminus \Sp_{\rm alg}(X))(K)$ is contained in the union over $i=1,\dots, r$ of the images of the composition
$$T_i(k) \cap U_i(K) \stackrel{\rho_i^\circ}{\lra}  Z_i(K)\lra Z(K)
\lra X(K).$$
Here the intersection $T_i(k) \cap U_i(K)$ is taken in 
$(T_{i,K})(K)$ via the canonical injection $T_i(k) \to (T_{i,K})(K)$.
\end{enumerate}
Moreover, if $X_{\ol K}$ does not contain any (possible singular) rational curve, then we can take every birational $K$-map 
$\rho_i:T_{i,K}\dashrightarrow Z_i$ 
to be a $K$-morphism under which $T_{i,K}$ is $Z_i$-isomorphic to the normalization $Z_i'$ of $Z_i$. 
\end{conj}

The essential case of the conjecture is when the transcendence degree of $K/k$ is 1, which implies the general case by induction.

\subsubsection*{Examples}

If $X_{\ol K}$ does not contain any rational curve, then the rational maps $\rho_i:T_{i,K}\dashrightarrow Z_i$ in Conjecture \ref{main conj copy} are claimed to be morphisms. The following simple example implies that this fails in the general case.

\begin{example}\label{blow up}
Let $K$ be the function field of one variable over $k=\CC$.
Set $Y=T_{K}$ for a smooth hyperbolic surface $T$ over $k$, and set $X$ to be the blowing-up of $Y$ along a closed point of $Y$ which is transcendental over $k$.
Then the birational map $T_{K}\dashrightarrow X$ is not a morphism.
\end{example}

In the following, we present an example to show that in Conjecture \ref{main conj copy}
(and Theorem \ref{partial to BL0 copy}), the finite set $\{Z_1,\dots, Z_r\}$ may contain more elements than the irreducible components of $Z$ to cover all but finitely many  rational points.

\begin{example}\label{exam isotrivial}
Let $C_0$ be a connected, smooth, and projective curve of genus greater than $1$ over $k=\CC$.
Let $K=k(C_0)$ be the function field, and denote by $C=C_0\times_kK$ the base change.

The finiteness of non-constant endomorphisms of $C$ implies that $\Sigma=C(K)\setminus C_0(k)$ is finite.  
%By the geometric Mordell conjecture proved by Manin \cite{Manin63}, the set
%$\Sigma=C(K)\setminus C_0(k)$ is finite.  
Moreover, $\Sigma$ is non-empty, since it contains the point $\Spec K\to C_0\times_kK$ induced by the canonical morphism $\Spec K\to C_0$ to the generic point of $C_0$.
Let $X=C\times_K C$ and $X_0=C_0\times_k C_0$ be the products. 
Then $X$ is algebraically hyperbolic and constant over $K$. 
The set $X(K)=C(K)\times C(K)$ is the union of the following four sets:
$$C_0(k)\times C_0(k),\qquad
C_0(k)\times  \Sigma, \qquad 
\Sigma\times C_0(k), \qquad 
\Sigma\times\Sigma.
$$   
Then $X(K)$ contains infinitely many points which are not from $X_0(k)$.
Moreover, for each $P\in \Sigma$, the subvariety $P \times_K C\simeq C$ of $X$ is abstractly isomorphic to $T_1(P)\times_kK$ for $T_1(P)=C_0$, so it is a constant variety.  
Similarly, 
$C \times_K P\simeq C$ is abstractly isomorphic to $T_2(P)\times_kK$ for $T_2(P)=C_0$.  
By these relations, $X(K)\setminus (\Sigma\times\Sigma)$ is covered by the $k$-points of the following varieties:
$$X_0,\qquad
T_1(P), \qquad 
T_2(P), \qquad
P\in \Sigma.
$$ 
\end{example}

\subsubsection*{When rational points are Zariski dense}

Here we present a weaker version of Conjecture \ref{main conj copy}. 
Conversely, this weaker version (for varying $X$) also implies Conjecture \ref{main conj copy}, which is similar to the fact that Theorem \ref{partial to BL00} implies Theorem \ref{partial to BL0 copy}.

\begin{conj} \label{main conj first step}
Let $K$ be a finitely generated field over a field $k$ of characteristic 0 such that $k$ is algebraically closed in $K$. 
Let $X$ be a projective variety over $K$. 
Assume that 
$(X\setminus \Sp_{\rm alg}(X))(K)$ is Zariski dense in $X$.
Then there is a birational $K$-map 
$\rho:T_{K}\dashrightarrow X$, 
where $T_{K}=T\times_k {K}$ is the base change for a projective variety $T$ over $k$.

Denote by $U$ a non-empty open $K$-subvariety of $T_{K}$ such that $\rho$ extends to a $K$-morphism $\rho^\circ:U\to X$.
Then the complement of $\Im(T(k) \cap U(K) \stackrel{\rho^\circ}{\to}  X(K))$
in $X(K)$ is not Zariski dense in $X$.
Here the intersection $T(k) \cap U(K)$ is taken in 
$(T_{K})(K)$ via the canonical injection $T(k) \to (T_{K})(K)$.

Moreover, if $X_{\ol K}$ does not contain any (possible singular) rational curve, then we can take the birational $K$-map 
$\rho:T_{K}\dashrightarrow X$ 
to be a $K$-morphism under which $T_{K}$ is $X$-isomorphic to the normalization $X'$ of $X$. 
\end{conj}

If we do not state where the rational points come from, then Conjecture \ref{main conj copy} has the following clean consequence (by combining the Green--Griffiths--Lang conjecture). 

\begin{conj} \label{main conj first step}
Let $K$ be a finitely generated field over a field $k$ of characteristic 0 such that $k$ is algebraically closed in $K$. 
Let $X$ be a projective variety over $K$. 
Assume that $X$ is of general type, and that $X$ is not birational constant over $K$.
Then $X(K)$ is not Zariski dense in $X$.
\end{conj}

Here a \emph{constant} variety over $K$ is the base change of a projective variety from $k$ to $K$. 
A \emph{birationally constant} variety over $K$ is a variety over $K$ which is birational to a constant variety over $K$.

\subsubsection*{Enhanced special set}

The spirit of Conjecture \ref{main conj copy} is that most rational points come from either abelian varieties or constant varieties. 
The contribution of abelian varieties is encoded in the special subset $\Sp_{\rm alg}(X)$, while the following conjecture asserts that we can enhance $\Sp_{\rm alg}(X)$ to contain the contribution of the constant varieties.

\begin{conj} \label{enhanced special set}
Let $K$ be a finitely generated field over an algebraically closed field $k$ of characteristic 0. 
Let $X$ be a projective variety over $K$. 
Then there is a unique Zariski closed subset $\Sp_{\rm alg+const}(X)$ of $X$ satisfying the following properties:
\begin{enumerate}[(1)]
\item $\Sp_{\rm alg+const}(X)$ contains $\Sp_{\rm alg}(X)$;
\item Every positive-dimensional birationally constant closed subvariety of $X$ is contained in $\Sp_{\rm alg+const}(X)$;
\item Every irreducible component of $\Sp_{\rm alg+const}(X)$ is positive-dimensional, and it is either birationally constant or an irreducible component of $\Sp_{\rm alg}(X)$.
\end{enumerate}
\end{conj}

In other words, $\Sp_{\rm alg+const}(X)$ is the unique maximal closed subset of $X$ whose irreducible components are either contained in $\Sp_{\rm alg}(X)$ or positive-dimensional birationally constant projective varieties over $K$.
Alternatively, $\Sp_{\rm alg+const}(X)$ is the union of $\Sp_{\rm alg}(X)$ with all positive-dimensional birationally constant closed subvarieties of $X$.  

It is interesting that Conjecture \ref{main conj copy} implies Conjecture \ref{enhanced special set} by taking $\Sp_{\rm alg+const}(X)$ to be the union of 
$\Sp_{\rm alg}(X)$ with all positive-dimensional irreducible components of the Zariski closure of $(X\setminus \Sp_{\rm alg}(X))(K)$ in $X$.

Based on the definition in Conjecture \ref{enhanced special set}, we see that 
Conjecture \ref{main conj copy} also implies the following conjecture.

\begin{conj} \label{main conj variant}
Let $K$ be a finitely generated field over an algebraically closed field $k$ of characteristic 0. 
Let $X$ be a projective variety over $K$. 
Then the set $(X\setminus \Sp_{\rm alg+const}(X))(K)$ is finite.
\end{conj}

\section{Canonical partial heights on abelian varieties}

The goal of this section is to use Betti maps to introduce canonical partial heights on abelian varieties, and prove the non-degeneracy of the canonical partial heights.

\subsection{Betti maps and Betti forms} \label{sec Betti}

Let us first review and set notations for Betti maps and Betti forms. 
These objects were introduced by Mok \cite[p. 374]{Mok1989} 
to study Mordell--Weil groups of abelian varieties over complex function fields. 
Our exposition mostly follows from \cite[\S2]{Cantat2021}, and we refer to the loc. cit. for more details. 

\subsubsection*{Betti maps}

Let $B$ be a Riemann surface.
Let $\pi:\CA\to B$ be a holomorphic family of abelian varieties over $B$; that is, $\CA$ is a complex manifold, $\pi$ is a smooth holomorphic map endowed with a holomorphic section $e:B\to\CA$, and every fiber of $\pi$ is an abelian variety with the identity point induced by $e$. 
Denote by $\CA(B)_h$ the group of holomorphic sections $s:B\to \CA$. 

If $B, \CA, \pi, e$ are algebraic, we will denote by 
$\CA(B)$ the group of algebraic sections $s:B\to \CA$. 
Then we have a natural injection $\CA(B)\to \CA(B)_h$.
However, in this subsection, we take the more general analytic setting.

Let $U\subset B$ be a connected and \emph{simply connected} open subset of $b$ in $B$. 
Let $b\in B$ be a point.
Denote by $\pi_U:\CA_U\to U$ the base change of $\pi$, and by $\CA_b$ the fiber of $\CA$ above $b$.
The \emph{Betti map} is a canonical real analytic map 
$$
\beta=\beta_{b,U}: \CA_U \lra \CA_b,
$$
which satisfies the following properties:
\begin{enumerate}[(1)]
\item The composition 
$\CA_{b}\hookrightarrow\CA_U\stackrel{\beta}{\to} \CA_b$ is the identity map. 
\item For any point $b'\in U$, the composition 
$\CA_{b'}\hookrightarrow\CA_U\stackrel{\beta}{\to} \CA_b$ is an isomorphism of real Lie groups.
\item The induced map 
$$\tilde\beta=(\beta,\pi): \CA_U \lra \CA_b\times U$$
 is a real analytic diffeomorphism of manifolds.
\item For any $x\in \CA_U$, the fiber $\beta^{-1}(\beta(x))$ is a complex analytic subset of $\CA_U$, which is biholomorphic to $U$. 
\end{enumerate}

The construction of the Betti map is as follows. 
The local system $R^1\pi_{*}\underline\ZZ$ is trivial on $U$ by the assumption that $U$ is simply connected. As a consequence, there are canonical isomorphisms $H^1(\CA_{b'},\ZZ)\to H^1(\CA_b,\ZZ)$ and 
$H_1(\CA_{b'},\ZZ)\to H_1(\CA_b,\ZZ)$ for any $b'\in U$. 
They induce real analytic isomorphisms
$$
\CA_{b'} \lra H_1(\CA_{b'},\RR)/H_1(\CA_{b'},\ZZ)
 \lra H_1(\CA_{b},\RR)/H_1(\CA_{b},\ZZ)
 \lra \CA_b.
$$
This gives the map in (2), which actually determine $\beta$.

The fiber $\CF_{x,U}=\beta_{b,U}^{-1}(\beta_{b,U}(x))$ is called \emph{the local Betti leaf at $x$ over $U$}. 
It is independent of the choice of $b\in U$ (for fixed $x,U$), and its 
germ $\CF_x$ at $x$ is independent of the choice of $U$. 
A \emph{Betti leaf} of $\CA$ is a connected subset $\CF_0$ of $\CA$ such that for any $x\in \CF_0$, the germ of $\CF_0$ at $x$ is equal to $\CF_x$. 
Note that any connected component of a torsion multi-section of $\CA\to B$ is a Betti leaf. 
The set of all Betti leaves forms a \emph{Betti foliation}. 

If $\CA$ is a trivial family over $B$ in that $\CA=\CA_0\times B$ for a complex abelian variety, then the Betti map is just the first projection from $\CA_U=\CA_0\times U$ to $\CA_b\simeq \CA_0$.
Then the maximal Betti leaves are just the horizontal sections, i.e. the fibers of the projection $\CA_0\times B \to \CA_0$.

\subsubsection*{Betti forms}

Let $\CL$ be a line bundle on $\CA$.  
For any $b\in B$, there is a unique translation-invariant closed $(1,1)$-form $\omega_b$ on $\CA_b$ representing the Chern class $c_1(\CL_b)$ on $\CA_b$. 
The \emph{Betti form $\omega=\omega(\CL)$ on $\CA$} associated to $\CL$ is the unique $(1,1)$-form $\omega$ on $\CA$ satisfying the following properties:
\begin{enumerate}[(1)]
\item For any $b\in B$, the pull-back of $\omega$ via $\CA_b\to \CA$ is equal to the translation-invariant K\"ahler form $\omega_b$. 
\item The pull-back of $\omega$ via the identity section $e:B\to \CA$
 is 0.
\end{enumerate}
These two properties determines $\omega$ uniquely. 
For the existence, let $(U,b)$ be a pair as above, and define the \emph{Betti form on $\CA_U$} associated to $\CL$ to be the differential form 
$\omega_U=(\beta,\pi)^*\omega_b$. 
It turns out that $\omega_U$ is a closed $(1,1)$-form on 
$\CA_U$. 
Moreover, $\omega_U$ is independent of the choice of $b\in U$ for fixed $U$, so we can glue $\omega_U$ for different $U$ to form the 
{Betti form $\omega$ on $\CA$}.

The Betti form $\omega$ is a closed $(1,1)$-form on $\CA$ satisfying the invariance property $[m]^*\omega=m^2\omega$. 
Because of this, one can also construct it by Tate's limiting argument.
Moreover, $\omega$ is semipositive if $\CL$ is ample on fibers of $\pi:\CA\to B$.

\subsection{Canonical partial heights} \label{sec partial}

The goal of this subsection is to introduce partial heights in terms of the Betti form. We will start with a review of the case of the N\'eron--Tate height.

Let $K=\CC(B)$ be the function field of a smooth quasi-projective curve $B$ over $\CC$. 
Let $A$ be an abelian variety over $K$ and let $L$ be a symmetric and ample line bundle over $A$. 
Denote by 
$$\hat h_L: A(\overline K)\lra \RR$$ 
the N\'eron--Tate height function defined by Tate's limiting argument. 

Replacing $B$ by an open subvariety if necessary, we assume that $A\to \Spec K$ extends to an abelian scheme $\pi:\CA\to B$. 
Note that we have assumed that $B$ is quasi-projective (instead of projective) in this section to get this flexibility.
We make an identification $\CA(B)= A(K)$ via the canonical isomorphism.  
In particular, 
we have the N\'eron--Tate height function  
$\hat h_L: \CA(B)\to \RR$.

Let $\CL$ be a line bundle on $\CA$ extending $L$. 
Let $\omega$ be the Betti form on $\CA$ associated to $\CL$. 
Note that $\omega$ is independent of the choice of $\CL$, since different choices of $\CL$ differ by a vertical divisor of $\CA$ and give isomorphic $\CL|_{\CA_b}$ for fixed $b\in B$.

We have the following interpretation of the canonical height in terms of Betti forms proved by
Gauthier--Vigny \cite[Thm. B]{GV2019}.

\begin{thm}\label{canonical height}
For any section $s:B\to \CA$,  
$$
\hat h_L(s)=\int_{B} s^*\omega. 
$$
\end{thm}

\subsubsection*{Canonical partial height}

Let $\omega$ be the Betti form on $\CA$ depending on a symmetric and ample line bundle $L$ as above. 
Let $D\subset B$ be a measurable subset.
%The measurability makes sense by taking a coordinate chart.
We will eventually take $D$ to be a disc under some local coordinate.
Recall that {the partial height of a section $s\in \CA(B)$ with respect to $(D,\omega)$} is
$$
h_{(D,\omega)}(s):=\int_{D} s^*\omega.
$$
To emphasize the use of the Betti form, we may also call $h_{(D,\omega)}(s)$ a \emph{canonical partial height of $s$}.

This gives a {canonical partial height function} $h_{(D,\omega)}: \CA(B)\to \RR$, and thus a canonical partial height function $h_{(D,\omega)}: A(K)\to \RR$.

If $D=B$, then we recover the original N\'eron--Tate height by Theorem \ref{canonical height}. 
In general, we always have 
$$
0\leq h_{(D,\omega)}(s)\leq h_{(\omega, B)}(s)=\hat h_L(s). 
$$

Similar to the N\'eron--Tate height, the canonical partial height is \emph{quadratic} in the following sense.

\begin{pro}[quadraticity]\label{quadratic}
\begin{enumerate}[(1)]
\item
For any points $s\in \CA(B)$ and $m\in \ZZ$,
$$
h_{(D,\omega)}([m]s)=m^2 h_{(D,\omega)}(s),
$$
\item
For any points $s,t\in \CA(B)$,  
$$
h_{(D,\omega)}(s+t)+h_{(D,\omega)}(s-t)=2h_{(D,\omega)}(s)+2 h_{(D,\omega)}(t).
$$
\end{enumerate}
\end{pro}

\begin{proof}
We refer to \cite[\S3]{se1} for the quadraticity of the classical N\'eron--Tate height $\hat h_L$. 
The proof for the canonical partial height is similar. 

We first prove (1). 
By the compatibility of the Betti maps and the Betti forms via the homomorphism $[m]:\CA\to \CA$, 
the symmetry condition $[m]^*L=L^{\otimes m^2}$ implies $[m]^*\omega=m^2\omega$.
The pull-back via $s:B\to \CA$ gives 
$$
(ms)^*\omega=m^2(s^*\omega). 
$$
Integrating over $D$, we obtain (1). 

For (2), start with the theorem of the cube, which asserts that
$$
m_{123}^*L \otimes (m_{12}^*L)^\vee\otimes (m_{23}^*L)^\vee \otimes (m_{13}^*L)^\vee\otimes m_{1}^*L \otimes m_{2}^*L \otimes m_{3}^*L=0
$$
in $\Pic(A^3)$. 
Here for any subset $I\subset \{1,2,3\}$, we have
$$
m_I:A^3\lra A, \quad (x_1,x_2,x_3)\longmapsto \sum_{i\in I}x_i. 
$$ 
By the compatibility of the Betti maps and the Betti forms, the relation implies 
$$
m_{123}^*\omega - m_{12}^*\omega- m_{23}^*\omega- m_{13}^*\omega+ m_{1}^*\omega+m_{2}^*\omega+m_{3}^*\omega=0
$$
as a differential form on $\CA\times_B \CA\times_B \CA$. 
For sections $s,t,u\in \CA(B)$, the pull-back via $(s,t,u):B\to \CA\times_B \CA\times_B \CA$ of the relation becomes
$$
(s+t+u)^*\omega - (s+t)^*\omega- (t+u)^*\omega- (s+u)^*\omega+ s^*\omega+t^*\omega+u^*\omega=0
$$
Setting $u=-t=[-1]\circ t$ in the above relation, we obtain 
$$
(s+t)^*\omega + (s-t)^*\omega= 2(s^*\omega+t^*\omega).
$$
Integrating over $D$, we obtain (2). 
\end{proof}

By the proposition, the canonical partial height defines a positive semi-definite quadratic form on $\CA(B)$. 
This induces a bilinear pairing on $\CA(B)$ by
$$
\pair{s,t}_{(D,\omega)}:=\frac12 \big( h_{(D,\omega)}(s+t) -h_{(D,\omega)}(s) -h_{(D,\omega)}(t)\big), \quad
s,t\in \CA(B). 
$$ 
By linearity, this extends to a bilinear pairing
$$
\pair{\cdot,\cdot}_{(D,\omega)}:\CA(B)_\RR\times \CA(B)_\RR \lra \RR
$$ 
and a quadratic form 
$$
h_{(D,\omega)}:\CA(B)_\RR \lra \RR, \quad 
s\longmapsto \pair{s,s}_{(D,\omega)}. 
$$ 
Based on the Betti form, we have the following result.

\begin{lem} \label{generalized form}
For any $s,t\in \CA(B)_\RR$, there is a unique $(1,1)$-form $\omega(s,t)$ on $B$ such that for any measurable subset $D$ of $B$, 
$$
\pair{s,t}_{(D,\omega)} = \int_D \omega(s,t). 
$$
The $(1,1)$-form satisfies the following extra properties:
\begin{enumerate}[(1)]
\item For any $s\in \CA(B)$, $\omega(s,s)=s^*\omega$. 
\item For any $s\in \CA(B)_\RR$, $\omega(s,s)$ is semipositive on $B$. 
\item The assignment $(s,t)\mapsto \omega(s,t)$ is $\RR$-bilinear in $s,t\in \CA(B)_\RR$. 
\item For any $s,t\in \CA(B)_\RR$,
$\omega(s,t)$ is real analytic on $B$ in the sense that 
under any complex analytic local coordinate $z$ on $B$, $\alpha=i f(z)dz\wedge d\bar z$ for a real analytic function $f(z)=f(\Re z, \Im z)$ (defined locally). 
\end{enumerate}

\end{lem}

\begin{proof}
The uniqueness of $\omega(s,t)$ is obtained by varying $D$ as discs in $B$. 
Now we prove the existence. 
If $s,t\in \CA(B)$, by the quadraticity in Proposition \ref{quadratic}, we can simply take 
$$
\omega(s,t)=\frac12 \big((s+t)^*\omega-s^*\omega-t^*\omega\big).
$$
If $s,t\in A(K)_\RR$ are general elements, 
write $s=a_1s_1+\dots + a_r s_r$ and 
$t=b_1t_1+\dots + b_{r'} t_{r'}$ for $s_i,t_j\in A(K)$ and
$a_i,\dots, b_j \in \RR$. 
Then we set 
$$
\omega(s,t)
=\sum_{i,j} a_i b_j \omega(s_i, s_j). 
$$
This satisfies the requirement by the quadraticity in Proposition \ref{quadratic},

For the extra properties, (1) and (3) follow from the uniqueness, and (4) follows from the construction since $\omega$ is real analytic. 
For (2), it holds for $s\in A(K)$ by (1) and the semipositivity of $\omega$, 
 it holds for $s\in A(K)_\QQ$ by passing to a multiple of $s$ in $A(K)$, 
and thus it holds for $s\in A(K)_\RR$ by approximation. 
This finishes the proof. 
\end{proof}

\subsubsection*{Non-degeneracy 1} 

Resume the above notations for $(B,K,A, L, \CA,\CL, \omega)$.
In particular, we have assumed that $L$ is symmetric and ample.

Denote by $(A^{(K/\CC)},\tr)$ Chow's $(K/\CC)$-trace of $A$. 
Then $A^{(K/\CC)}$ is an abelian variety over $\CC$ together with a $K$-homomorphism 
$\tr: A^{(K/\CC)}\times_\CC K\to A$, and $A^{(K/\CC)}$ is the unique final object with such property.
Since we are in characteristic 0, the trace map $\tr: A^{(K/\CC)}\times_\CC K\to A$ is actually a closed immersion.
We refer to \cite{Conrad} for more details on Chow's trace. 
The following Lang--N\'eron theorem is a combination of \cite[Thm. 2.1, Lem. 9.13, Thm. 9.15]{Conrad}.

\begin{thm}[Lang--N\'eron] \label{LN}
The abelian group $A(K)/A^{(K/\CC)}(\CC)$ is finitely generated.
The N\'eron--Tate height 
function $\hat h_L: A(K) \to \RR$ is invariant under the translation action of $A^{(K/\CC)}(\CC)$ on $A(K)$, and induces a positive definite quadratic form
$$
\hat h_L: (A(K)/A^{(K/\CC)}(\CC))\otimes_\ZZ\RR \lra \RR.
$$
\end{thm}

The main result of this section is the following non-degeneracy theorem, which is a ``partial'' version of the Lang--N\'eron theorem. 

\begin{thm}[non-degeneracy] \label{non-deg}
Let $D$ be an open disc in $B$. 
Then the canonical partial height function $ h_{(D,\omega)}: A(K) \to \RR$ is invariant under the translation action of $A^{(K/\CC)}(\CC)$ on $A(K)$, and induces a positive definite quadratic form
$$
h_{(D,\omega)}: (A(K)/A^{(K/\CC)}(\CC))\otimes_\ZZ\RR \lra \RR.
$$
As a consequence, there is a constant $\epsilon>0$ such that 
$$
h_{(D,\omega)}(s)\geq \epsilon \hat h_{L}(s), \quad \forall s\in A(K).
$$
\end{thm}

\begin{proof}

The last statement follows from the abstract situation that both $h_{(D,\omega)}$ and $\hat h_{L}$ are positive definite quadratic forms on the finite-dimensional vector space 
$(A(K)/A^{(K/\CC)}(\CC))\otimes_\ZZ\RR$.

Now we prove that $\hat h_{(D,\omega)}: A(K) \to \RR$ is invariant under the translation action of $A^{(K/\CC)}(\CC)$ on $A(K)$. 
The key is that $\hat h_{(D,\omega)}(s)=0$ for any $s\in A^{(K/\CC)}(\CC)$. 
This is a consequence of the inequality $\hat h_{(D,\omega)}(s)\leq \hat h_L(s)$ and Theorem \ref{LN}; alternatively, 
by functoriality, we can reduce it to the case $A=A^{(K/\CC)}\times_\CC K$, and in this case the Betti form $\omega$ is the pull-back of a differential form from $A^{(K/\CC)}$.

Once $\hat h_{(D,\omega)}$ is zero on $A^{(K/\CC)}(\CC)$, it is invariant under the translation action of $A^{(K/\CC)}(\CC)$ on $A(K)$ by some abstract properties of quadratic forms. 
In fact, since $\omega$ is semipositive, the values of $\hat h_{(D,\omega)}: A(K) \to \RR$ are non-negative, so $\hat h_{(D,\omega)}: A(K)_\RR \to \RR$ is positive semi-definite (on any finite-dimensional subspace).
Then the pairing $\pair{\cdot,\cdot}_{(D,\omega)}$ satisfies the Cauchy--Schwartz inequality 
$\pair{s,t}_{(D,\omega)}^2 \leq \pair{s,s}_{(D,\omega)} \pair{t,t}_{(D,\omega)}$
for $s,t\in A(K)$. 
Then $\pair{s,t}_{(D,\omega)}=0$ for any $s\in A^{(K/\CC)}(\CC)$ and any $t\in A(K)$, and thus $\hat h_{(D,\omega)}(s+t)=\hat h_{(D,\omega)}(t)=0$ by the quadraticity. 
This proves that $\hat h_{(D,\omega)}: A(K) \to \RR$ is invariant under the translation action of $A^{(K/\CC)}(\CC)$ on $A(K)$.

It remains to prove that $h_{(D,\omega)}$ is positive definite on
$(A(K)/A^{(K/\CC)}(\CC))\otimes_\ZZ\RR$.
Let $s\in A(K)_\RR$ be an element with $h_{(D,\omega)}(s)=0$.
We need to prove $s\in A^{(K/\CC)}(\CC)_\RR$.
By Lemma \ref{generalized form}, there is a semi-positive and real analytic $(1,1)$-form $\alpha=\omega(s,s)$ on $B$ such that
$$
h_{(D,\omega)}(s)=\int_D \alpha.
$$
As $\alpha$ is semipositive, $\alpha|_D=0$ on $D$. 
As $\alpha$ is real analytic, $\alpha|_D=0$ on $D$ forces $\alpha=0$ on $B$. 
In fact, a basic result asserts that the zero locus of a nonzero real analytic function on a connected real analytic manifold has measure 0.
We refer to \cite{Mityagin} for this basic fact for connected open sets of $\RR^n$, which can be extended to connected real analytic manifolds by taking a cover by open balls.

By Lemma \ref{generalized form} again, 
$$\hat h_L(s)=h_{(\omega, B)}(s)=\int_B \alpha=0.$$ 
Then $s\in A^{(K/\CC)}(\CC)_\RR$ by Theorem \ref{LN}. 
This finishes the proof. 
\end{proof}

\subsubsection*{Non-degeneracy 2} 

Theorem \ref{non-deg} confirms the non-degeneracy conjecture (Conjecture \ref{conj non-deg}) for abelian varieties, and implies the conjecture for varieties finite over abelian varieties.
For convenience, we state and prove the result in the following.

\begin{thm}  \label{confirm conj non-deg}
Let $K=\CC(B)$ be the function field of a smooth projective curve $B$ over $\CC$. 
Let $X$ be a projective variety over $K$  with a finite morphism $f:X\to A$ for an abelian variety $A$ over $K$. 
Then Conjecture \ref{conj non-deg} holds for $X$. 
\end{thm}

\begin{proof}

Resume the notations $(B,K, X, L, h_L, D, \omega, \{x_n\})$ of the conjecture. Note that the curve $B$ is assumed to be projective. 

We first confirm the conjecture for the case that $X$ is an abelian variety.
Denote by $U$ a non-empty open subvariety of $B$ over which $\CX$ is smooth.
By Theorem \ref{non-deg}, the conjecture holds for the pair $(D, \omega)$, where $D$ is an open disc in $U$, and $\omega$ is a Betti form on $\CX_U$.
By Lemma \ref{comparison partial}(1), the conjecture holds for any strictly positive pair $(D, \omega)$ on $\CX$ such that $D$ is an open disc in $U$.
The restriction that $D$ is an open disc contained in $U$ is easy to remove. 
In fact, we can also shrink $D$ to make it contained in $U$, and this process decreases the partial height. 
This confirms the conjecture for abelian varieties.

For the general case $f:X\to A$, 
let $\CA\to B$ be an integral model of $A$ over $B$.
The morphism $f:X\to A$ extends to a morphism $\tilde f: \CA_U\to \CX_U$
for some non-empty open subvariety $U$ of $B$.
To prove the conjecture, we can assume $L=f^*L_0$ for an ample line bundle $L_0$ on $A$. 
Let $(D_0,\omega_0)$ be a strictly positive pair on $\CA$, where $D_0$ is an open disc in $U$. 
We have 
$$h_{L}(x_n)=h_{L_0}(f(x_n))+O(1) $$
and 
$$
h_{(D_0,\tilde f^*\omega_0)}(x_n)=h_{(D_0,\omega_0)}(f(x_n)).
$$
Then the conjecture holds for $(D_0,\tilde f^*\omega_0)$. 
We can adjust this to any strictly positive pair 
$(D,\omega)$ on $\CX$, as in the case of abelian varieties.  
\end{proof}

\section{Hyperbolic covers of abelian varieties}

The goal of this section is to review some preliminary results on special sets and then prove Theorem \ref{hyperbolic main}. 
The strategy is a variant of the proof of Theorem \ref{partial to BL0 copy}, and the key is the non-degeneracy of partial heights in Theorem \ref{confirm conj non-deg}.

\subsection{Some preliminary results} \label{sec preliminary}

In this section, we first review some preliminary results on hyperbolicity and generic sequences.

\subsubsection*{Equality of special sets} 

The hyperbolicity theory for varieties finite over abelian varieties is very complete due to the work of Ueno, Kawamata, Noguchi, Winkelmann and Yamanoi. 
%The hyperbolicity theory for varieties finite over abelian varieties is very complete due to the work of Ueno, Kawamata and Yamanoi. 
We summarize the relevant result as follows.

\begin{thm} \label{special set}
Let $X$ be a projective variety over $\CC$ which has a finite morphism to an abelian variety. 
Then the following hold:
\begin{enumerate}[(1)]
\item $\Sp_{\rm alg}(X)=\Sp_{\rm an}(X)$.   
\item $X$ is of general type if and only if $\Sp_{\rm an}(X)\neq X$.  
\end{enumerate}
\end{thm}
\begin{proof}
We first see that if $X$ is not of general type, then
$\Sp_{\rm alg}(X)=\Sp_{\rm an}(X)= X$. 
This is a direct consequence of the structure theorem in \cite[Thm. 13]{Kawamata81}, which asserts that the normalization of $X$ has an \'etale cover of the form $X_0\times A_0$ for a projective variety $X_0$ and an abelian variety $A_0$ of positive dimension. 

This proves the ``if'' part of (2). 
The ``only if'' part of (2) follows from Yamanoi \cite[Cor. 1]{Yamanoi15} (see \cite{Noguchi2007} for some previous results).
The loc. cit. assumes that $X$ is smooth, but it is extended to the singular case by Hironaka's resolution of singularities.

Now we prove (1). 
By definition, $\Sp_{\rm alg}(X)\subseteq\Sp_{\rm an}(X)$.
To prove the inverse inclusion, by definition, 
$\Sp_{\rm an}(Y)=Y$ for every irreducible component $Y$ of $\Sp_{\rm an}(X)$. 
Then (2) implies that $Y$ is not of general type, so 
$\Sp_{\rm alg}(Y)=Y$ by (2). It follows that $Y = \Sp_{\rm alg}(Y)\subseteq \Sp_{\rm alg}(X)$, and the union of $Y$ gives 
This proves that $\Sp_{\rm an}(X) \subseteq \Sp_{\rm alg}(X)$. 
\end{proof}

The following basic result asserts that taking the special subset of a projective variety is stable under base change, see  \cite[Proposition 3.7]{Javanpeykar2022} for a variant.

\begin{lem}\label{special set base change}
Let $K'$ be a field extension of a field $K$.
Let $X$ be a projective variety over $K$. 
Then $\Sp_{\rm alg}(X_{K'})=\Sp_{\rm alg}(X)_{K'}$ as closed subsets of $X_{K'}$.
\end{lem}
\begin{proof}
By definition, $\Sp_{\rm alg}(X)_{K'}\subseteq \Sp_{\rm alg}(X_{K'})$. We need to prove the opposite direction. 
Let $G$ be an abelian variety over a finite extension $K''$ of $K'$ together with a non-constant rational map $\lambda:G\dashrightarrow X_{K''}$. 
We need to prove that the image of $G$ in $X_{K'}$ is contained in $\Sp_{\rm alg}(X)_{K'}$, or equivalently, 
the image of $G$ in $X_{K''}$ is contained in $\Sp_{\rm alg}(X)_{K''}$. 
By the Lefschetz principle, we can assume that $K''$ 
is finitely generated extension of $K$.

Let $S$ be a variety over $K$ with function field $K''$. 
Replacing $S$ by an open subvariety if necessary, we can assume that $G\to \Spec K''$ extends to an abelian scheme $\CG\to S$, and 
$\lambda:G\dashrightarrow X_{K''}$ extends to a rational map 
$\bar\lambda:\CG \dashrightarrow X\times_K S$ whose indeterminacy locus does not contain any fiber of $\CG\to S$. 
By definition, for any closed point $s\in S$, the image of the composition 
$\CG_s\to \CG \dashrightarrow X\times_K S \to X$
is contained in $\Sp_{\rm alg}(X)$. 
Note that the Zariski closure of $\cup_s \CG_s$ in $\CG$ is $\CG$. 
It follows that $\Sp_{\rm alg}(X)$ actually contains the image of 
the composition 
$\CG \dashrightarrow X\times_K S \to X$. 
In particular, $\Sp_{\rm alg}(X)$ contains the image of 
the composition 
$$G\lra \CG \stackrel{\bar \lambda}{\dashrightarrow} X\times_K S \lra X,$$
which is isomorphic to the composition 
$$G\stackrel{\lambda}{\dashrightarrow}  X_{K''} \lra  X.$$
After base change, $\Sp_{\rm alg}(X)_{K''}$ contains the image of 
$\lambda:G\dashrightarrow X_{K''}$. 
This finishes the proof.
\end{proof}

Now we have the following consequence of the above results.

\begin{cor} \label{special set2}
Let $X$ be a projective variety over a field $K$ of characteristic 0 which has a finite morphism to an abelian variety over $K$. 
Then $X$ is of general type if and only if $\Sp_{\rm alg}(X)\neq X$.  
\end{cor}
\begin{proof}
Let $f:X\to A$ be  a finite morphism to an abelian variety $A$ over $K$. 
By the Lefschetz principle, we can descend $(X,A,f)$ to a finitely generated subfield $K_0$ of $K$ and take a base change by an embedding $K_0\to \CC$. 
Then the result follows from Theorem \ref{special set} and Lemma \ref{special set base change}.
\end{proof}

\subsubsection*{Transcendental closed points}

We first introduce the notion of transcendental closed points on complex curves, which have ``functorial'' specialization properties.

Let $B$ be a smooth curve over $\CC$, and denote by $K=\CC(B)$ the function field.
By the Lefschetz principle, there is a finitely generated subfield $k_0$ of $\CC$ and a smooth curve $B_0$ over $k_0$ such that $B\simeq B_0\times_{k}\CC$. 
Denote $K_0=k_0(B_0)$.

A point $b\in B(\CC)$ is called a \emph{transcendental closed point with respect to $B_0/k_0$} if the image of the composition
$$
\Spec\CC \stackrel{b}{\lra} B\lra   B_0 
$$
is the generic point of $B_0$.  
Note that we have canonical bijections
$$
B(\CC)\simeq \Hom_{\Spec \CC}(\Spec \CC, B)
\simeq \Hom_{\Spec k_0}(\Spec \CC, B_0).
$$
Under the bijections, $b\in B(\CC)$ is transcendental if it is given by a composition 
$$
\Spec\CC \stackrel{\lambda_b}{\lra} \Spec K_0 \stackrel{\eta_0}{\lra} 
B_0,
$$
where $\eta_0:\Spec K_0\to B_0$ is the generic point, and $\lambda_b:\Spec\CC \to \Spec K_0$ is induced by an embedding $K_0\to \CC$ of fields over $k_0$.
The basic result is that ``almost all'' points of $B(\CC)$ are transcendental.

\begin{lem} \label{transcendental0}
 There are only countably many points $b\in B(\CC)$ which are not {transcendental with respect to $B_0/k_0$}.
\end{lem}
\proof
Let $b\in B(\CC)$ be a point which is not transcendental. 
Then the image of the corresponding map 
$\Spec \CC\to B_0$ is a closed point of $B_0$. Each such closed point gives rise to finitely many $b$, and there are only countably many closed points of $B_0$. 
This finishes the proof.
\endproof

As above, let $B$ be a smooth curve over $\CC$, and denote by $K=\CC(B)$ the function field.
Now let $\CX$ be a variety over $\CC$ with a flat projective morphism $\pi:\CX\to B$ over $\CC$.
Denote by $X=\CX_K$ the generic fiber of $\CX\to B$.
We will see that for ``almost all'' $b\in B(\CC)$, the fiber $\CX_b$ is  related to $X$ in a convenient functorial manner. 

Now we choose the model $B_0/k_0$ to accommodate $\CX$. Namely, 
let $k_0$ be a finitely generated subfield of $\CC$ such that $(B, \CX, \pi)$ is the base change of a triple $(B_0, \CX_0, \pi_0)$ from $k_0$ to $\CC$. 
Denote $K_0=k_0(B_0)$ and $X_0=\CX_{0,K_0}$.

Let $b\in B(\CC)$ be a {transcendental closed point with respect to $B_0/k_0$} given by the composition 
$$
\Spec\CC \stackrel{\lambda_b}{\lra} \Spec K_0 \stackrel{\eta_0}{\lra} 
B_0. 
$$
Then the fiber of $\CX$ above $b$ is given by 
$$
\CX_b=\CX\times_{B} (\Spec\CC,b)
\simeq (\CX_0\times_{B_0} \Spec K_0)\times_{\Spec K_0} (\Spec\CC,\lambda_b),$$
so we have a canonical isomorphism
$$
\CX_b \simeq X_0 \times_{\Spec K_0} (\Spec\CC,\lambda_b).
$$
In summary, the construction from $X$ to $\CX_b$ can be re-interpreted as the process $X \to X_0\leftarrow \CX_b$, i.e., a descent from $K$ to $K_0$ followed by a base change from $K_0$ to $\CC$. 
This process brings many good properties, and the following is one of them. 

\begin{lem} \label{transcendental1}
\begin{enumerate}[(1)]
\item For any  {transcendental closed point $b\in B(\CC)$ with respect to $B_0/k_0$}, the special set $\Sp_{\rm alg}(\CX_b)$ is isomorphic to some base change of some descent of $\Sp_{\rm alg}(X)$. 
\item Denote by $\overline{\Sp_{\rm alg}(X)}$ the Zariski closure of $\Sp_{\rm alg}(X)$ in $\CX$. Assume furthermore that $\overline{\Sp_{\rm alg}(X)}\hookrightarrow \CX$ can also be descended to $B_0$.
Then for any  {transcendental closed point $b\in B(\CC)$ with respect to $B_0/k_0$}, $\Sp_{\rm alg}(\CX_b)$ is equal to the fiber of $\overline{\Sp_{\rm alg}(X)}$ above $b$.
\end{enumerate}
\end{lem}
\proof
These are consequences of the construction and Lemma \ref{special set base change}.
\endproof

The concept of transcendental closed points applies to finitely many varieties and morphisms. Namely, let $B$ be a smooth curve over $\CC$, and $\CX_i$ for $i=1,\dots, r$ 
be a variety over $\CC$ with a morphism $\pi_i:\CX_i\to B$ over $\CC$, and let $f_j:\CX_{a(j)}\to \CX_{a'(j)}$ for $j=1,\dots, m$ and $a(j),a'(j)\in \{1,\dots, r\}$ be a $B$-morphism among them. 
By the Lefschetz principle, there is a finitely generated subfield $k_0$ of $\CC$ such that $(B, \pi_i:\CX_i\to B, f_j:\CX_{a(j)}\to \CX_{a'(j)})$ is the base change of $(B_0, \pi_{i,0}:\CX_{i,0}\to B_0, f_{j,0}:\CX_{a(j),0}\to \CX_{a'(j),0})$ from $k_0$ to $\CC$. 
Then for any transcendental closed points of $b\in B(\CC)$ with respect to $B_0/k_0$, the specializations
$\CX_{i,b}$ and $f_{j,b}:\CX_{a(j),b}\to \CX_{a'(j),b}$
are base changes of descents of 
$X_{i}$ and $f_{j,K}:X_{a(j)}\to X_{a'(j)}$. 
Here $K=\CC(B)$ and $X_i=\CX_{i,K}$ as usual.

\subsection{Rational points} \label{sec hyperbolic}

In this subsection, we prove Theorem \ref{hyperbolic main}. For convenience, it is duplicated in the following with a detailed statement. 

\begin{thm}[Theorem \ref{hyperbolic main}] \label{hyperbolic main copy}
Let $K$ be a finitely generated field over a field $k$ of characteristic 0 such that $k$ is algebraically closed in $K$. 
Let $X$ be a projective variety over $K$ with a finite morphism $f:X\to A$ for an abelian variety $A$ over $K$. 
Assume that $X$ is algebraically hyperbolic in that $\Sp_{\rm alg}(X)=\emptyset$. 
Let $Z$ be the Zariski closure of 
$X(K)$ in $X$. 

Then Conjecture \ref{main conj} holds for $X/K/k$.
Namely, let $Z$ be the Zariski closure of 
$X(K)$ in $X$. 
Then there is a finite set $\{Z_1,\dots, Z_r\}$ of distinct closed subvarieties of $Z$
containing all irreducible components of $Z$ and satisfying the following conditions: 
\begin{enumerate}[(1)]
\item
For each $i$, the normalization $Z_i'$ of $Z_i$ is constant in the sense that
there is a $K$-isomorphism 
$\rho_i:T_i\times_k {K}\to Z_i'$ for a projective variety $T_i$ over $k$.
\item 
The set $X(K)$ is equal to the union over $i=1,\dots, r$ of the images of the composition
$$T_i(k)\lra (T_i\times_k {K})(K)\stackrel{\rho_i}{\lra} Z_i'(K) \lra Z_i(K)\lra Z(K)
\lra X(K).$$
\end{enumerate}
\end{thm}

The theorem is essentially (but not directly) implied by Theorem \ref{partial to BL0 copy}, without some extra arguments to convert the different setting. First, Theorem \ref{hyperbolic main copy} is still a consequence of its special case $Z=X$ as follows.

\begin{thm} \label{hyperbolic copy}
Let $K$ be a finitely generated field over a field $k$ of characteristic 0 such that $k$ is algebraically closed in $K$. 
Let $X$ be a projective variety over $K$ with a finite morphism $f:X\to A$ for an abelian variety $A$ over $K$. 
Assume that $X$ is algebraically hyperbolic in that $\Sp_{\rm alg}(X)=\emptyset$. 
Assume that $X(K)$ is Zariski dense in $X$. 
Then the normalization of $X$ is isomorphic to 
the base change $T_K=T\times_k {K}$ for a normal projective variety $T$ over $k$, and the complement in $X(K)$ of the image of the composition
$T(k)\to (T_{K})(K)
\to X(K)$
is not Zariski dense in $X$. 
\end{thm}

The implication of Theorem \ref{hyperbolic main copy} by Theorem \ref{hyperbolic copy} is the same as the implication of Theorem \ref{partial to BL0 copy} by Theorem \ref{partial to BL00}. So we omit it here.

\subsubsection*{Reduce the transcendence degree} 

The remaining part of this subsection is devoted to prove Theorem \ref{hyperbolic copy}. 
Since the case $K=k$ is trivial, we assume that $\trdeg(K/k)\geq 1.$
We first reduce  
$\trdeg(K/k)$ to 1 in the theorem. 
The key is the following lemma,  a well-known consequence of Chow's $K/k$-trace. 

\begin{lem} \label{trace argument}
Let $K/k$ be an extension fields of characteristic 0 such that $k$ is algebraically closed in $K$. Let 
$f:T\times_k K \to A$
be a $K$-morphism, where $T$ is a projective variety over $k$ and $A$ is an abelian variety over $K$. 
Assume that there is a $k$-point $t_0\in T(k)$ contained in the regular locus of $T$.
Then the image of $f$ is contained in the translation $(A^{K/k})_K+a$ of $(A^{K/k})_K$ by $a=f(t_0)\in A(K)$. 
Moreover, the composition
$$T\times_k K \lra (A^{K/k})_K+a \stackrel{-a}{\lra} (A^{K/k})_K$$ 
is equal to the base change of a $k$-morphism 
$T\to A^{K/k}$. 
\end{lem}

\begin{proof}
We first reduce to the case that $T$ is smooth over $k$.
Let $T'\to T$ be a resolution of singularities, which exists by Hironaka's theorem. 
We can further assume that $T'\to T$ is an isomorphism above the regular point $t_0\in T(k)$, and thus $t_0$ has a unique preimage $t_0'\in T'(k)$.
Denote by $f':T'_{K}\to A$ the composition of
$T'_{K} \to T_{K}\to A$.
Replacing $(T, t_0, f)$ by $(T', t_0', f')$, we can assume that $T$ is smooth over $k$.

Now we assume that $T$ is smooth over $k$.
Replacing $f:T_K\to A$ by $f-f(t_0)$, 
we can assume that $a=f(t_0)=0$.
Let $\mathrm{Alb}( T)$ be the Albanese variety of $( T,t_0)$ over $k$, endowed with the $k$-morphism $ T\to \mathrm{Alb}( T)$ sending $t_0$ to 0.
Then $f:T_{K}\to A$ factorizes through a $K$-homomorphism $\mathrm{Alb}( T)_K\to A$.
The latter further factorizes through the base change of a $k$-homomorphism $\mathrm{Alb}( T)\to A^{K/k}$.
This gives the result.
\end{proof}

Now it is easy to reduce Theorem \ref{hyperbolic copy} to the case $\trdeg(K/k)=1$.
In fact, assume that $\trdeg(K/k)>1$, and we need to prove that the conclusion holds for the datum $(K/k, X, A)$. 
Assume that $X$ is normal by taking a normalization.

Let $k_1$ be an intermediate field of $K/k$ such that $K/k_1$ has transcendence degree 1 and that $k_1$ is algebraically closed in $K$.
Assume that the theorem holds for $(K/k_1, X, A)$. 
Then there is a $K$-isomorphism 
$T_1\times_{k_1} {K}\to X$ for a projective variety $T_1$ over $k_1$, and the complement of $\Im(T_1(k_1)\to X(K))$ in $X(K)$ is not Zariski dense in $X$. 

By Lemma \ref{trace argument}, the finite $K$-morphism 
$T_1\times_{k_1} K\to A$ induces a finite $k_1$-morphism 
$T_1\to A^{K/k_1}$. 
This reduces the problem to $(k_1/k, T_1, A^{K/k_1})$. 
Repeat the process on $k_1/k$, we eventually get the case $\trdeg(K/k)=1$.
This finishes the reduction process.

In the following, we assume $\trdeg(K/k)=1$ in the proof of Theorem \ref{hyperbolic copy}. 
Then the case of $k=\CC$ of the theorem is a special case of Theorem \ref{partial to BL00}.
For general $k$, the statement that the normalization of $X$ is isomorphic to 
the base change $T_K$ can be deduced from Theorem \ref{partial to BL00}
by the Lefschetz principle and a descent argument. 
However, the statement about the image of 
$T(k)\to X(K)$
does not seem a direct consequence of Theorem \ref{partial to BL00}.
Due to this, we will sketch a full path of the proof of Theorem \ref{hyperbolic copy}
(assuming $\trdeg(K/k)=1$).

%The idea is similar to the proof of Theorem \ref{partial to BL00}. 

\subsubsection*{Dominated by a constant family} 

We want to apply the proof of Theorem \ref{partial to BL00} to conclude that $X$ is dominated by a constant family, but the problem is that the base field $k$ is more general than $\CC$.
To remedy the problem, we will alter the proof to the  current setting.

Let $X/A/K/k$ be as in Theorem \ref{hyperbolic copy}, and assume $\trdeg(K/k)=1$. 
Let $h_L:X(\ol K)\to \RR$ be a Weil height function with respect to an ample line bundle $L$ on $X$. 
We claim that $h_L$ is bounded  on $X(K)$. 

If $k=\CC$, the claim is obtained in the proof of Theorem \ref{partial to BL00}.
For that, we need to fulfill the conditions of Theorem \ref{partial to BL00}. 
First, Conjecture \ref{conj non-deg} for the non-degeneracy of partial heights holds for the current $X$ by Theorem \ref{confirm conj non-deg}.
Second, for any integral model $\CX\to B$ of $X$ over $B$, there is an open disc $D$ in $B$ such that all fibers of $\CX$ above $D$ are hyperbolic. 
In fact, by assumption, $\Sp_{\rm alg}(X)=\emptyset$, so 
 Lemma \ref{transcendental0} and Lemma \ref{transcendental1} imply that there is a point $b\in B(\CC)$ with 
$\Sp_{\rm alg}(\CX_b)=\emptyset$. 
By Theorem \ref{special set}, $\Sp_{\rm an}(\CX_b)=\emptyset$. 
By a consequence of the Brody lemma (cf. \cite[III, Prop. 3.1]{Lang87}), being  hyperbolic is an open condition under the Euclidean topology.
Then we can find a disc $D$ with center $b$ satisfying the requirement.

For general $k$, we can easily deduce the claim by the Lefschetz principle. 
In fact, it suffices to prove that for any infinite sequence $\{x_n\}_{n\geq1}$ of $X(K)$, the sequence $\{h_L(x_n)\}_{n\geq1}$ is  bounded. 
The datum $(K,X,A, f:X\to A)$ is defined over a finitely generated subfield of $k$ over $\QQ$, and so is every $x_n\in X(K)$.  
Then all these data with all $n\geq1$ are defined over a countably generated subfield $k_0$ of $k$ over $\QQ$. 
Fix an inclusion $k_0\hookrightarrow \CC$. 
By descent from $k$ to $k_0$ and base change from $k_0$ to $\CC$, we achieve the case that the base field is $\CC$.
The process does not change the heights of points, so we conclude that $\{h_L(x_n)\}_{n\geq1}$ is bounded. 
This proves the claim in the general case.

Once the height is bounded, the argument applying Hilbert schemes holds over $k$ (instead of $\CC$), so the first step in the proof of Theorem \ref{partial to BL00} holds in the current situation. 
Namely, there is a disjoint union $V$ of finitely many smooth projective varieties  over $k$ together with a surjective $K$-morphism $V\times_k {K} \to X$ such that the  composition
$$V(k)\lra (V\times_k {K})(K)
\lra X(K)$$
is surjective. 
Note that in the current situation, the morphism $V\times_k {K} \to X$ is also extended from a rational map $V_K\dashrightarrow X$. This is possible still because $X$ does not contain any rational curve, as a consequence of the finite morphism $f:X
\to A$.

\subsubsection*{Constancy} 

Here we continue to prove that the normalization of $X$ is isomorphic to $T_K$ for some projective variety $T$ over $k$.
We can assume that $X$ is normal, since the general case can be solved by a normalization process.

As mentioned above, this can be proved by Theorem \ref{constancy} by applying the Lefschetz principle. Note that Theorem \ref{constancy} is proved by a rather involved foliation method. 
However, in the current situation, taking advantage of the finite morphism $f:X\to A$, we have an algebraic proof in this case, which is significantly easier than the original proof. 
Due to this merit, we will provide it in the following.

Replacing $V$ by a connected component of $V$, we assume that $V$ is connected. 
Thus we have a smooth projective variety $V$ over $k$ together with a surjective $K$-morphism $\varphi:V_{K} \to X$ such that $\varphi(V(k))$ is Zariski dense in $X$.
Consider the composition
$$V_K\stackrel{\varphi}{\lra}  X \stackrel{f}{\lra} A.$$
Take a point $v_0\in V(k)$, which exists since $\varphi(V(k))$ is Zariski dense in $X$.
By Lemma \ref{trace argument}, 
up to replacing $f:X\to A$ by a translation by an element of $A(K)$, 
the morphism $V_K\to A$ factorizes through the base change of a $k$-morphism $V\to A^{K/k}$. Denote by $W$ the image of the finite morphism $V\to A^{K/k}$. 
Then $W_K$ (as a closed subvariety of $A$) is the image of $V_K\to A$, and also  the image of $X\to A$.
In summary, we have a composition 
$$
V_K \stackrel{\varphi}{\lra} X \lra W_K \longrightarrow (A^{K/k})_K \lra A
$$  
of $K$-morphisms, where the first two arrows are surjective, 
the last three arrow are finite, 
and the composition
$V_K\to W_K$ is the base change of a $k$-morphism $V\to W$. 

Denote by $V'$ the normalization of $V\to W$. 
Then $V'_K$ is still normal, and thus the normalization of $V_K\to W_K$. 
Moreover, $V'_K$ still dominates $X$, since $X$ is finite over $W_K$. 
Replacing $V$ by $V'$, we can assume that $V_K\to X$ is finite (and surjective). Here we might lose the smoothness of $V$, but this will not be used below. 
Denote by $k(V)^c$ the Galois closure of $k(V)/k(W)$. 
Replacing $V$ by the normalization of $W$ in $k(V)^c$, we can assume that 
$k(V)/k(W)$ is Galois. 

By uniqueness of normalization,  the action of 
$\Gal(k(V)/k(W))$ on $k(V)$ extends to an action on $V/W$. 
Via the natural isomorphism 
$$\Gal(k(V)/k(W))\lra \Gal(K(V_K)/K(W_K)),$$ 
the action of $\Gal(K(V_K)/K(W_K))$ on $V_K$ is through the action of  
$\Gal(k(V)/k(W))$ on $V$. 
The subgroup $\Gal(K(V_K)/K(X))$ of $\Gal(K(V_K)/K(W_K))$ corresponds to a subgroup $G$ of $\Gal(k(V)/k(X))$, and the action of $\Gal(K(V_K)/K(X))$ on $V_K$ is through the action of $G$ on $V$. 
Taking quotient, we have a composition
$$
(V/G)\times_k K \lra V_K/\Gal(K(V_K)/K(X))\lra X.
$$
Here the first arrow is a canonical isomorphism, and the second arrow is a birational and finite morphism. 
Note that $V$ is normal, and the quotient varieties are normal by construction, so the second arrow is isomorphic to the normalization of $X$. 
This prove that $X\simeq T_K$ with $T=V/G$.

\subsubsection*{Constancy of the points}

It remains to prove that via $X\simeq T_{K}$, the complement of 
$\Im(T(k)\to X(K))$ in $X(K)$
is not Zariski dense in $X$. 
Similar to the proof of Theorem \ref{partial to BL00}, this is a consequence of 
Theorem \ref{unique constant structure}.

%%%%%%%%%
\bibliographystyle{alpha}
\bibliography{dd}

%%%%%%%%%
\medskip
\noindent \small{Address: \textit{Beijing International Center for Mathematical Research, Peking University, Beijing 100871, China}}

\noindent \small{Email: \textit{xiejunyi@bicmr.pku.edu.cn}}

\medskip
\noindent \small{Address: \textit{Beijing International Center for Mathematical Research, Peking University, Beijing 100871, China}}

\noindent \small{Email: \textit{yxy@bicmr.pku.edu.cn}}

%------------------------------------------------------
\end{document}